\documentclass[a4paper,UKenglish,cleveref]{lipics-v2021}
\pdfoutput=1 

\usepackage{mathtools} 
\usepackage{stmaryrd} 
\usepackage{tikz}
\usetikzlibrary{matrix,arrows,calc,backgrounds,shapes,shapes.geometric,patterns,positioning,fit,decorations.pathmorphing,decorations.pathreplacing,shadows}

\usepackage{todonotes}

\newcommand{\cutout}[1]{}

\theoremstyle{definition}
\newtheorem{Definition}[theorem]{Definition}

\renewcommand{\phi}{\varphi}
\renewcommand{\theta}{\vartheta}
\renewcommand{\emptyset}{\varnothing}
\renewcommand{\epsilon}{\varepsilon}

\newcommand{\FO}{{\rm FO}}

\newcommand*{\tup}[1]{\bar{#1}}
\newcommand{\ta}{\tup a}
\newcommand{\tb}{\tup b}
\newcommand{\tx}{\tup x}

\newcommand{\tz}{\tup z}

\renewcommand{\AA}{{\mathfrak A}}
\newcommand{\BB}{{\mathfrak B}}
\newcommand{\N}{{\mathbb N}} 
\newcommand{\B}{{\mathbb B}} 
\newcommand{\K}{{\mathbb K}} 

\DeclareMathOperator{\Lit}{\mathrm{Lit}}

\newcommand{\Bool}{\mathbb{B}}

\newcommand{\Nat}{\mathbb{N}}

\newcommand{\Trop}{\mathbb{T}}
\newcommand{\Vit}{\mathbb{V}}

\newcommand{\PosBool}{\mathsf{PosBool}}

\renewcommand{\bar}{\overline}

\newcommand{\absorb}{\succeq}

\tikzset{
    arr/.style={draw,->,>=stealth',shorten <=2pt,shorten >=2pt,every node/.style={auto,inner sep=2pt,font=\scriptsize}},
    gamenode/.style={draw,inner sep=2pt,minimum size=.4cm},
    p0/.style={gamenode,circle},
    p1/.style={gamenode,rectangle},
    F/.style={thick, pattern=north east lines},
    dot/.style={circle,draw,fill,black,minimum size=3pt,inner sep=0pt},
    marker/.style={draw=none,inner sep=0pt,overlay},
    short/.style={ shorten >=#1, shorten <=#1 }
}

\newcommand*{\StrA}{\mathfrak{A}}
\newcommand*{\StrB}{\mathfrak{B}}
\newcommand*{\Rat}{\mathbb{Q}}
\newcommand*{\Real}{\mathbb{R}}
\newcommand*{\Sorb}{\mathbb{S}}
\newcommand*{\Why}{\mathbb{W}}

\newcommand*{\from}{\colon}
\newcommand*{\isoto}{\xrightarrow{\sim}}
\newcommand*{\Lcomp}{\overline{L}}

\newcommand*{\ttov}{\sigma_{\Trop \to \Vit}}
\newcommand*{\vtot}{\sigma_{\Vit \to \Trop}}

\DeclareMathOperator{\id}{id}
\DeclareMathOperator{\Hom}{Hom}
\DeclareMathOperator{\arity}{arity}

\newcommand*{\iso}{\mathrel{\cong}}

\DeclarePairedDelimiter{\intpr}{\llbracket}{\rrbracket}

\newtheorem*{question*}{Question}
\newtheorem*{questions*}{Questions}

\newcommand*{\Fraisse}{Fraïssé}
\newcommand*{\EF}{Ehrenfeucht--\Fraisse{}}

\makeatletter
\tikzset{
        hatch distance/.store in=\hatchdistance,
        hatch distance=7pt,
        hatch thickness/.store in=\hatchthickness,
        hatch thickness=2pt
        }
\pgfdeclarepatternformonly[\hatchdistance,\hatchthickness]{north east hatch}
    {\pgfqpoint{-200pt}{-200pt}}
    {\pgfqpoint{\hatchdistance+199pt}{\hatchdistance+199pt}}
    {\pgfpoint{\hatchdistance-1pt}{\hatchdistance-1pt}}%
    {
        \pgfsetcolor{\tikz@pattern@color}
        \pgfsetlinewidth{\hatchthickness}
        \pgfpathmoveto{\pgfqpoint{0pt}{0pt}}
        \pgfpathlineto{\pgfqpoint{\hatchdistance}{\hatchdistance}}
        \pgfusepath{stroke}
    }
\pgfdeclarepatternformonly[\hatchdistance,\hatchthickness]{north west hatch}
    {\pgfqpoint{-200pt}{-200pt}}
    {\pgfqpoint{\hatchdistance+199pt}{\hatchdistance+199pt}}
    {\pgfpoint{\hatchdistance-1pt}{\hatchdistance-1pt}}%
    {
        \pgfsetcolor{\tikz@pattern@color}
        \pgfsetlinewidth{\hatchthickness}
        \pgfpathmoveto{\pgfqpoint{\hatchdistance}{0pt}}
        \pgfpathlineto{\pgfqpoint{0pt}{\hatchdistance}}
        \pgfusepath{stroke}
    }
\makeatother

\title{Elementary equivalence versus isomorphism in semiring semantics }

\author{Erich Grädel}{RWTH Aachen University, Germany}{graedel@logic.rwth-aachen.de}{}{}
\author{Lovro Mrkonjić}{RWTH Aachen University, Germany}{mrkonjic@logic.rwth-aachen.de}{}{}

\authorrunning{E. Grädel and L. Mrkonjić}

\Copyright{Erich Grädel, Lovro Mrkonjić}



\nolinenumbers 
\hideLIPIcs  

\begin{document}

\maketitle

\begin{abstract}
We study the first-order axiomatisability of finite semiring interpretations or, equivalently, the
question whether elementary equivalence and isomorphism coincide for valuations
of atomic facts over a finite universe into a commutative semiring. Contrary to the classical case
of Boolean semantics, where every finite structure can obviously be axiomatised up to isomorphism
by a first-order sentence, the situation in semiring semantics is rather different, and strongly depends on
the underlying semiring. We prove that for a number of important semirings, including min-max semirings,
and the semirings of positive Boolean expressions, there exist finite semiring interpretations that
are elementarily equivalent but not isomorphic. The same is true
for the polynomial semirings $\Sorb[X]$, $\Bool[X]$ and $\Why[X]$ that are
universal for the classes of absorptive, idempotent, and fully idempotent semirings, respectively.
On the other side, we prove that for other, practically relevant, semirings such as the Viterby semiring $\Vit$,
the tropical semiring $\Trop$, the natural semiring $\N$ and the universal polynomial
semiring $\N[X]$, all finite semiring interpretations are first-order axiomatisable (and thus elementary equivalence implies isomorphism),
although some of the axiomatisations that we exhibit use an infinite set of axioms.
\end{abstract}

\section{Introduction}

Semiring semantics is based on the idea to evaluate logical statements not just by \emph{true}
or \emph{false}, but by values in some commutative semiring $(K,+,\cdot,0,1)$.
In this context, the standard semantics appears as the special case 
when the Boolean semiring $\Bool = (\{\bot, \top\}, \lor, \land, \bot, \top)$ is used.
Valuations in other semirings provide additional information, beyond the truth
or falsity of a statement:
the Viterbi-semiring $\Vit = ([0,1]_{\Real}, \max, \cdot, 0, 1)$ models \emph{confidence scores}, 
the tropical semiring $\Trop= (\Real_{+}^{\infty}, \min, +, \infty, 0)$
is used for \emph{cost analysis}, and min-max-semirings $(K, \max, \min, a, b)$ for a 
totally ordered set $(K,<)$ can model, for instance, different \emph{access levels}.
More generally, semirings of polynomials, such as $\N[X]$, $\Why[X]$ or $\Bool[X]$, 
allow us to track the role of specific atomic facts for the evaluation of 
a logical statement, to describe evaluation strategies for a formula, and
to determine which combinations of literals prove the truth of a formula.

Some of the motivation for the study of semiring semantics comes from the successful development
of semiring provenance in database theory and related fields (see e.g. \cite{DeutchMilRoyTan14, DoleschalKimMarPet20,GeertsPog10,
GreenKarTan07, GreenTan17, OzakiPen18, RaghothamanMenZhaNaiSch20, Senellart17}),
and the fact that the typical applications of provenance analysis, such as confidence scores, cost analysis,
proof counting, and the understanding of evaluation strategies are of importance in many other
areas of logic as well.
However, semiring provenance analysis for database queries had originally been largely
confined to positive query languages, such as conjunctive queries, positive relational algebra, and Datalog,
and the treatment of negation poses non-trivial algebraic problems. Only recently, semiring
semantics has been extended to logics with negation, and in particular to full first-order logic
\cite{GraedelTan17, GraedelTan20}, by means of new algebraic constructions based on quotient semirings.
Semiring semantics has also been studied for other logics, including modal logics, description logics,
guarded logics, and fixed-point logic \cite{BourgauxOzaPenPre20, DannertGra19, DannertGra20, DannertGraNaaTan21},
and this paper is part of a larger project devoted to a systematic study of semiring semantics for various logics. 
An important objective in this context is the understanding of the \emph{model theory of semiring semantics},
and the development of model-theoretic methods for semiring interpretations.

It turns out that this is much more involved and diverse than for Boolean semantics.
In the standard semantics, a model $\AA$ 
assigns to each (instantiated) atomic formula a Boolean value,
whereas $K$-interpretations $\pi$, for a suitable semiring $K$, generalise this
by assigning to each literal a semiring value in $K$,
where $0$ is interpreted as \emph{false} and all other semiring values as \emph{nuances of true}.
Interpreting disjunction by $+$ and conjunction by $\cdot$, we can extend $\pi$ to
provide semiring valuations $\pi \intpr{\phi}\in K$ for all first-order sentences $\phi$, written
in negation normal form.
Semiring semantics thus 
gives a finer distinction of logical statements, and formulae that are equivalent in the Boolean sense
(i.e. in the Boolean semiring) may have different valuations in other semirings. As a consequence,
standard facts of classical (finite) model theory may lead to interesting and sometimes rather difficult questions
in semiring semantics, and the answer may strongly depend on algebraic properties of the underlying
semiring. Specific such questions that we study here concern the first-order axiomatisability of finite
$K$-interpretations or, what amounts to the same, the relationship between isomorphism and 
elementary equivalence in this context. 

It is a rather trivial fact of finite model theory that every finite 
structure $\AA$ (with a finite vocabulary $\tau$) can be axiomatised, 
up to isomorphism, by a first-order sentence $\chi_\AA$. 
In particular, two finite $\tau$-structures $\AA$ and $\BB$ are isomorphic if,
and only if, they are elementarily equivalent, in short $\AA \equiv \BB$, 
which means that they cannot be 
distinguished by any first-order sentence. Is this also the case for semiring interpretations?
Notice that standard notions such as isomorphism and elementary equivalence generalise in
a natural way from $\tau$-structures to semiring interpretations, which raises,
for any given semiring $K$, the following 

\begin{questions*}
\begin{enumerate}
\item Are elementary equivalent finite $K$-interpretations always isomorphic?\label{item:question1}
\item Is every finite $K$-interpretation $\pi_A$ first-order axiomatisable, in the
sense that there is a set of axioms $\Phi_A \subseteq \FO$ such that whenever
$\pi_B \intpr{\phi} = \pi_A \intpr{\phi}$ for all $\phi \in \Phi_A$, then $\pi_B \cong \pi_A$?\label{item:question2}
\item Does every finite $K$-interpretation admit an axiomatisation by a \emph{finite} set of axioms?\label{item:question3}
\item Can every finite $K$-interpretation be axiomatised by a single first-order sentence?\label{item:question4}
\end{enumerate}
\end{questions*}

Clearly, the first two questions are equivalent, and a
positive answer to the third question implies also positive ones to the first two.
The converse is not necessarily true, because a first-order axiomatisation of a finite
semiring interpretation might require an infinite collection of sentences,
and, contrary to the Boolean case,
it is a priori also not clear that an axiomatisation by a finite set of sentences implies
an axiomatisation by a single sentence, because from the
value of a conjunction we cannot necessary infer the values of its components.

We shall prove that the answers to these questions strongly depend on the chosen semiring.
There are in fact rather simple semirings, such as min-max semirings with at least three elements,
for which one can construct examples of non-isomorphic $K$-interpretations which
are, however, elementarily equivalent. The standard method for proving elementary equivalence
in model theory, the \EF{} method, seems not really available in
semiring semantics, an aspect that we shall discuss at the end of this paper.
To establish elementary equivalence, we shall hence develop
new methods based on classes of semiring homomorphisms and reduction arguments.
Elementarily equivalent but non-isomorphic semiring interpretations also exist
for powerful polynomial semirings, such as $\Sorb[X]$, $\Bool[X]$ and $\Why[X]$ which are
universal for the classes of absorptive, idempotent, and fully idempotent semirings, respectively.
On the other side, there are practically relevant semirings, such as the Viterby semring $\Vit$,
the tropical semiring $\Trop$, the natural semiring $\N$ and the universal polynomial
semiring $\N[X]$, for which any finite $K$-interpretation is first-order axiomatisable,
thus elementary equivalence does indeed imply
isomorphism. At least for $\Vit$ and $\Trop$, finite axiomatisations are always possible,
but not axiomatisations by a single sentence, so there exist semirings where the answers
to \hyperref[item:question3]{questions (3) and (4)} are different.

\section{Semiring interpretations}

We briefly summarise semiring semantics for first-order logic, as introduced in \cite{GraedelTan17}.

\begin{Definition}[Semiring]
A commutative semiring $\K = (K, +, \cdot, 0, 1)$ is an algebraic structure with two binary operations
such that $(K, +, 0)$ and $(K, \cdot, 1)$ are commutative monoids,
multiplication distributes over addition and multiplication with zero annihilates elements.
We refer to commutative semirings as semirings and identify $\K$
with its universe $K$ if the operations are clear from the context.
\end{Definition}

Let $\tau$ denote a finite relational vocabulary. We write $\Lit_n(\tau)$ for the set of atoms $R\tz$ 
and negated atoms $\neg R\tz$ with $R \in \tau$ and where $\tz$ is any tuple of variables taken
from $\{x_1, \dots, x_n\}$. For a universe $A$, we write $\Lit_A(\tau)$ for the set of
\emph{instantiated} $\tau$-literals $R\ta$ and $\neg R\ta$ with $\ta \in A^{\arity(R)}$.

 \begin{Definition}[$K$-interpretation]
For a semiring $K$, a mapping $\pi \from \Lit_A(\tau) \to K$ is called a $K$-\emph{interpretation}
over the universe $A$ with signature $\tau$. We call $\pi$ \emph{model-defining} if
exactly one of the values $\pi(L)$ and $\pi(\Lcomp)$ is zero
for all pairs of opposing literals $L, \Lcomp \in \Lit_A(\tau)$.
In that case, $\pi$ \emph{induces} the classical model $\StrA_{\pi}$
with $\StrA_{\pi} \models L$ if, and only if, $\pi(L) \neq 0$.
\end{Definition}

\begin{Definition}[Isomorphism]
$K$-interpretations $\pi_A \from \Lit_A(\tau) \to K$ and $\pi_B \from \Lit_B(\tau) \to K$ are \emph{isomorphic},
denoted as $\pi_A \iso \pi_B$, if there is a bijective mapping $\sigma \from A \to B$ such that
\[\pi_A(L) = \pi_B(\sigma(L)) \quad \text{for all } L \in \Lit_A(\tau),\]
where $\sigma(L) \in \Lit_B(\tau)$ is defined by replacing each $a \in A$ occurring in $L$ with $\sigma(a) \in B$.
The mapping $\sigma$ is called an \emph{isomorphism}, denoted as $\sigma \from \pi_A \isoto \pi_B$.
\end{Definition}

Given a $K$-interpretation $\pi \from \Lit_A(\tau) \to K$, a formula $\phi(\tx) \in \FO(\tau)$ in negation normal form
and an assignment $\ta \subseteq A$, the semiring semantics $\pi \intpr{\phi(\ta)}$ is straightforwardly defined
by induction on $\FO(\tau)$. We first extend $\pi$ by mapping equalities and inequalities to their truth values by 
\[\pi \intpr{a = b} \coloneqq\begin{cases} 1 &\text{ if } a = b \\0 &\text{ if } a \neq b \end{cases}\quad\text{and}\quad
 \pi \intpr{a \neq b} \coloneqq\begin{cases} 0 &\text{ if } a = b \\1 &\text{ if } a \neq b \end{cases},\]
and by interpreting disjunctions and existential quantifiers as sums, and conjunctions and universal quantifiers as products:
\begin{alignat*}{3}
\pi \intpr{\psi(\ta) \lor \theta(\ta)} &\coloneqq \pi \intpr{\psi(\ta)} + \pi \intpr{\theta(\ta)} &\quad\quad\quad \pi \intpr{\psi(\ta) \land \theta(\ta)} &\coloneqq \pi \intpr{\psi(\ta)} \cdot \pi \intpr{\theta(\ta)} \\
\pi \intpr{\exists x \theta(\ta, x)} &\coloneqq \sum_{a \in A} \pi \intpr{\theta(\ta, a)} &\quad\quad\quad \pi \intpr{\forall x \theta(\ta, x)} &\coloneqq \prod_{a \in A} \pi \intpr{\theta(\ta, a)}.
\end{alignat*}
Our goal is to analyse model-theoretic concepts from classical model theory under semiring semantics.

\begin{Definition}[Elementary Equivalence]
\phantomsection\label{def:equivalence}
Let $\pi_A \from \Lit_A(\tau) \to K$ and $\pi_B \from \Lit_B(\tau) \to K$ be two $K$-interpretations,
$\ta \in A^k$ and $\tb \in B^k$ be $k$-tuples and $k \in \Nat$.
The pairs $\pi_A, \ta$ and $\pi_B, \tb$ are \emph{elementarily equivalent}, denoted as $\pi_A, \ta \equiv \pi_B, \tb$, if
\[\pi_A \intpr{\phi(\ta)} = \pi_B \intpr{\phi(\tb)} \quad\text{for all } \phi(\tx) \in \FO(\tau) \text{ with } k \text{ free variables.}\]
They are $m$-\emph{equivalent} with $m \in \Nat$, denoted as $\pi_A, \ta \equiv_m \pi_B, \tb$,
if the above holds for all suitable $\phi(\tx)$ with quantifier rank at most $m$.
\end{Definition}

Clearly, the notions of isomorphism and elementary equivalence of $K$-interpretations are natural generalisations
of the corresponding definitions for $\tau$-structures. Further, it is obvious that, as in classical semantics,
isomorphism implies elementary equivalence.

\begin{lemma}[Isomorphism Lemma]
\label{lem:isomorphism}
Let $\pi_A \from \Lit_A(\tau) \to K$ and $\pi_B \from \Lit_B(\tau) \to K$ be two $K$-interpretations,
$\ta \in A^k$ and $\tb \in B^k$ be $k$-tuples and $\sigma \from \pi_A \isoto \pi_B$ an isomorphism with $\sigma(\ta) = \tb$.
Then, $\pi_A \intpr{\phi(\ta)} = \pi_B \intpr{\phi(\tb)}$ holds for all $\phi(\tx) \in \FO(\tau)$ with $k$ free variables.
\end{lemma}

Coarser definitions may be conceivable in semirings $K$ with more than two elements,
such as replacing equality by a congruence relation ${\sim} \subseteq K \times K$.
However, \hyperref[def:equivalence]{Definition~\ref*{def:equivalence}} indirectly
covers these variants, since any non-trivial congruence relation $\sim$ on $K$
induces a semiring homomorphism $h_{\sim} \from K \to K/{\sim}$,
which is compatible with $\FO$-semantics as follows \cite{GraedelTan17}.

\begin{lemma}[Fundamental Property]
\label{lem:fundamental}
Let $\pi \from \Lit_A(\tau) \to K$ be a $K$-interpretation and $h \from K \to L$ a semiring homomorphism.
Then, $(h \circ \pi)$ is an $L$-interpretation and
\[(h \circ \pi) \intpr{\phi(\ta)} = h(\pi \intpr{\phi(\ta)})\]
holds for all $\phi(\tx) \in \FO(\tau)$ and $\ta \subseteq A$. Thus, the following diagram commutes.

\vspace{\baselineskip}
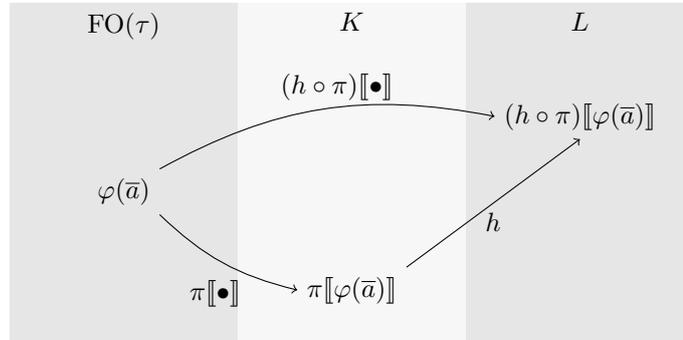
\begin{minipage}{\linewidth}
\centering
\begin{tikzpicture}
\fill[black!10] (0cm, -0.5cm) rectangle (3cm, -5cm);
\fill[black!3] (3cm, -0.5cm) rectangle (6cm, -5cm);
\fill[black!10] (6cm, -0.5cm) rectangle (9cm, -5cm);
\node[anchor=north](fo) at (1.5cm, -0.5cm) {$\FO(\tau)$};
\node[anchor=north](k) at (4.5cm, -0.5cm) {$K$};
\node[anchor=north](l) at (7.5cm, -0.5cm) {$L$};
\node[anchor=center](phi) at (1.5cm, -3cm) {$\phi(\ta)$};
\node[anchor=north](piphi) at (4.5cm, -4cm) {$\pi \intpr{\phi(\ta)}$};
\node[anchor=center](hpiphi) at(7.5cm, -2cm) {$(h \circ \pi) \intpr{\phi(\ta)}$};
\draw (phi.south east) edge[->, bend right=15] node[midway, below, xshift=-0.1cm, yshift=-0.1cm] {$\pi \intpr{\bullet}$} (piphi.west);
\draw(piphi.north east) edge[->] node[midway, below] {$h$} (hpiphi.south);
\draw(phi.north east) edge[->, bend left=20] node[midway, above, xshift=0.2cm] {$(h \circ \pi) \intpr{\bullet}$} (hpiphi.west);
\end{tikzpicture}
\end{minipage}
\end{lemma}

\section{Polynomial semirings and the universal property}

Semiring homomorphisms and \hyperref[lem:fundamental]{the fundamental property}
open the possibility to reduce semiring semantics to the evaluation of polynomials.
For a finite set $X$ of abstract provenance tokens that are used to track atomic facts,
consider the semiring $\Nat[X]$ of multivariate polynomials with indeterminates from $X$ and coefficients from $\Nat$,
whose generality is formalised by the following \emph{universal property} \cite{GreenKarTan07}.

\begin{lemma}[Universal Property]
\label{lem:universal}
For each commutative semiring $K$, every assignment $e \from X \to K$ induces
a unique homomorphism $h_e \from \Nat[X] \to K$ with $h_e(x) = e(x)$ for $x \in X$.
\end{lemma}

To demonstrate the use of this property, consider the $\Vit$-interpretation
$\pi_{\Vit}$ over $A = \{a, b\}$ depicted below on the left
and a corresponding $\Nat[X]$-interpretation $\pi$ on the right,
which maps all the true literals to their own variable in $X \coloneqq \{w, x, y, z\}$.

\vspace{\baselineskip}
\begin{minipage}{\linewidth}
\centering
$\pi_{\Vit}:\quad$
\begin{tabular}{c | c | c | c | c |}
$A$ & $P$ & $Q$ & $\neg P$ & $\neg Q$ \\ \hline
$a$ & $0.3$ & $0.9$ & $0$ & $0$ \\
$b$ & $0$ & $0.5$ & $0.4$ & $0$ \\
\end{tabular}
$\quad\quad\quad\pi:\quad$
\begin{tabular}{c | c | c | c | c |}
$A$ & $P$ & $Q$ & $\neg P$ & $\neg Q$ \\ \hline
$a$ & $x$ & $y$ & $0$ & $0$ \\
$b$ & $0$ & $z$ & $w$ & $0$ \\
\end{tabular}
\end{minipage}
\vspace{\baselineskip}

Given, for example, $\psi \coloneqq \forall x (Px \lor Qx)$, we get
$\pi_{\Vit} \intpr{\psi} = \max \{0.3, 0.9\} \cdot \max\{0, 0.5\} = 0.45$,
but if we have already computed the polynomial
$p \coloneqq \pi \intpr{\psi} = (x + y) \cdot (0 + z) = xz + yz$,
then plugging in the values $e \from x \mapsto 0.3, y \mapsto 0.9, z \mapsto 0.5$ and $w \mapsto 0.4$
induces the homomorphism $h_e \from \Nat[X] \to \Vit$, and evaluating $p$ gives the same value
$h_e(p) = \max\{0.3 \cdot 0.5, 0.9 \cdot 0.5\} = 0.45$
by \hyperref[lem:fundamental]{the fundamental property}.
This approach may be used to save computation resources
if we have several interpretations sharing the same set of true literals,
so that $\pi \intpr{\psi}$ is re-usable,
but the universality of $\Nat[X]$ also makes it relevant for model theory.

Other polynomial semirings can be used to capture smaller classes of semirings,
in the sense of \hyperref[lem:universal]{the universal property}.

\begin{Definition}[Idempotence and Absorption]
A semiring $K$ is called \emph{idempotent} if $a + a = a$ holds for all $a \in K$, that is, if addition is idempotent.
It is \emph{multiplicatively idempotent} if $a \cdot a = a$ for all $a \in K$. If both properties hold, we call $K$ \emph{fully idempotent}.
Finally, a semiring $K$ is \emph{absorptive} if $a + ab = a$ holds for all $a, b \in K$.
\end{Definition}

Note that absorptive semirings are idempotent and if we replace addition with $\lor$ and multiplication with $\land$,
absorption corresponds to the absorptive law in lattices.
In particular, absorptive and fully idempotent semirings are precisely the semirings induced by distributive lattices \cite{Naaf19}.

\begin{itemize}
\item By dropping coefficients from $\Nat[X]$, we get the semiring
$\B[X]$ whose elements are just finite sets of distinct monomials.
It has \hyperref[lem:universal]{the universal property} for the class of idempotent semirings.
\item By dropping also exponents, we get the semiring $\Why[X]$ of 
finite sums of monomials that are linear in each argument,
with \hyperref[lem:universal]{the universal property} for fully idempotent semirings.
It is sometimes called the Why-semiring.
\item For absorptive semirings, we require absorptive polynomials as introduced in \cite{DeutchMilRoyTan14}. 
An absorptive polynomial is a sum of distinct monomials over a finite set of variables $X$,
with absorption among monomials:
a monomial $m_1$ absorbs $m_2$, denoted $m_1 \absorb m_2$, if it has smaller exponents, 
i.e. if $m_1(x) \le m_2(x)$ for all $x \in X$, where $m(x)$ denotes the exponent of $x$ in $m$. 
Notice that absorption is the \emph{inverse} pointwise order on the exponents.
For example, $xy^2 \absorb x^3y^2$ and $x \absorb xy$, but $x^2y$ and $xy^2$ are incomparable.
In an absorptive polynomial, we omit all monomials that would be absorbed, so absorptive polynomials
are $\absorb$-\emph{antichains} of monomials (which are always finite \cite{GraedelTan20}). 
Consequently, addition and multiplication are defined as usual,
but afterwards we drop all monomials that are absorbed.
We write $\Sorb[X]$ for the semiring of absorptive polynomials over the finite variable set $X$
with the aforementioned operations. The $0$-element is the empty polynomial and
$1$ denotes the polynomial consisting just of the monomial $1$ (with all zero exponents).
It is not difficult to verify that this defines an absorptive semiring
which has \hyperref[lem:universal]{the universal property} for the class of all absorptive semirings. 
\item By dropping exponents from $\Sorb[X]$, we obtain yet another polynomial semiring 
$\PosBool[X]$, which is universal for the fully idempotent and absorptive semirings.
Incidentally, polynomials from $\PosBool[X]$, such as $x + yz$, represent all positive Boolean expressions
over the variables $X$ up to equivalence, if we replace addition by disjunction and multiplication by conjunction,
hence the name $\PosBool[X]$. This is the distributive lattice freely generated by the set $X$.
\end{itemize}

The following figure, adapted from \cite{Naaf19},
shows the relationships between the aforementioned classes of semirings,
together with their respective \hyperref[lem:universal]{universal polynomial semirings}.

\vspace{\baselineskip}
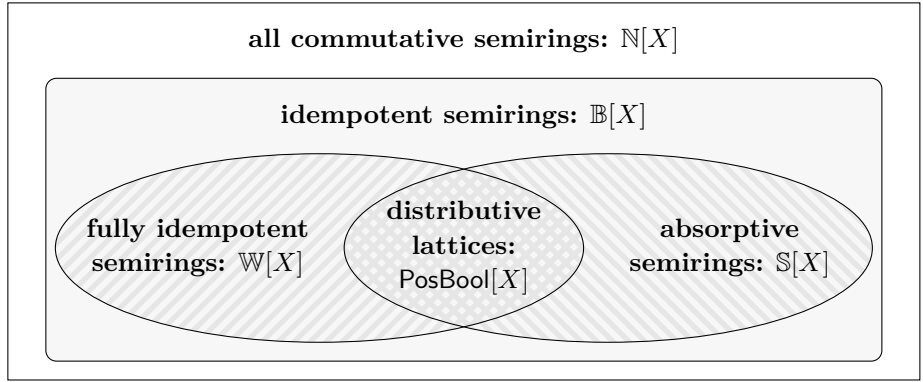
\begin{minipage}{\linewidth}
\centering
\begin{tikzpicture}
\draw (0cm, 0cm) rectangle (12cm, -5cm);
\fill[black!3, rounded corners] (0.5cm, -1cm) rectangle (11.5cm, -4.75cm);
\pattern[pattern=north east hatch, pattern color=black!10] (4.1cm, -3.25cm) ellipse (3.475cm and 1.25cm);
\pattern[pattern=north west hatch, pattern color=black!10] (7.9cm, -3.25cm) ellipse (3.475cm and 1.25cm);
\draw[rounded corners] (0.5cm, -1cm) rectangle (11.5cm, -4.75cm);
\draw (4.1cm, -3.25cm) ellipse (3.475cm and 1.25cm);
\draw (7.9cm, -3.25cm) ellipse (3.475cm and 1.25cm);
\node[anchor=center] at (6cm, -0.5cm) {\textbf{all commutative semirings:} $\Nat[X]$};
\node[anchor=center] at (6cm, -1.5cm) {\textbf{idempotent semirings:} $\Bool[X]$};
\node[anchor=center, align=center] at (6cm, -3.25cm) {\textbf{distributive} \\ \textbf{lattices:} \\ $\PosBool[X]$};
\node[anchor=center, align=center] at (2.5cm, -3.25cm) {\textbf{fully idempotent} \\ \textbf{semirings:} $\Why[X]$};
\node[anchor=center, align=center] at (9.5cm, -3.25cm) {\textbf{absorptive} \\ \textbf{semirings:} $\Sorb[X]$};
\end{tikzpicture}
\end{minipage}

\section{Separating elementary equivalence from isomorphism}
\label{sec:separation}

We shall now, for certain semirings $K$, provide examples of finite, non-isomorphic $K$-interpretations that are,
however, elementarily equivalent. We thus provide negative answers to \hyperref[item:question1]{Question (1)} from the introduction,
and hence also to \hyperref[item:question2]{Questions (2) to (4)}.
For instance, we claim that the following two $K_4$-interpretations
over the min-max-semiring with four elements, $K_4 = \{0, 1, 2, 3\}$, are elementarily equivalent, but not isomorphic.

\vspace{\baselineskip}
\begin{minipage}{\linewidth}
\centering
$\pi_{PQ}:\quad$
\begin{tabular}{c | c | c | c | c |}
$A$ & $P$ & $Q$ & $\neg P$ & $\neg Q$ \\ \hline
$a$ & 1 & 3 & 0 & 0 \\
$b$ & 2 & 1 & 0 & 0 \\
$c$ & 3 & 2 & 0 & 0 \\
\end{tabular}
$\quad\quad\quad\pi_{QP}:\quad$
\begin{tabular}{c | c | c | c | c |}
$A$ & $P$ & $Q$ & $\neg P$ & $\neg Q$ \\ \hline
$a$ & 3 & 1 & 0 & 0 \\
$b$ & 1 & 2 & 0 & 0 \\
$c$ & 2 & 3 & 0 & 0 \\
\end{tabular}\hspace*{8mm}
\end{minipage}
\vspace{\baselineskip}

Observe that $\pi_{QP}$ is obtained from $\pi_{PQ}$ by permuting the relations $P$ and $Q$,
or visually, by permuting the ``columns''. Moreover, for both interpretations,
the $Q$-column can be obtained by permuting the $P$-column.
Informally, these properties ensure that the two interpretations are ``sufficiently similar''
so that no first-order sentence can distinguish them.

Clearly, $\pi_{PQ}$ is not isomorphic to $\pi_{QP}$, as intended. However,
the only tool we presented so far for proving elementary equivalence under semiring semantics
is \hyperref[lem:isomorphism]{the Isomorphism Lemma} itself, which is not directly applicable for obvious reasons.
Hence, we shall develop another tool for proving elementary equivalence
that enables the indirect use of \hyperref[lem:isomorphism]{the Isomorphism Lemma}
after switching to a different semiring via homomorphisms.

\subsection{Separating pairs of homomorphisms}

The central idea of the reduction technique is to ``decompose'' the semiring $K$ via homomorphisms.
Observe that if $\pi_A \not\equiv \pi_B$, then there is a witnessing sentence $\psi \in \FO(\tau)$
with $\pi_A \intpr{\psi} \neq \pi_B \intpr{\psi}$, hence a pair of distinct elements $a, b \in K$ with
$\pi_A \intpr{\psi} \eqqcolon a \neq b \coloneqq \pi_B \intpr{\psi}$ exists.
If we can find two homomorphisms $h_A, h_B \from K \to L$ with $h_A(a) \neq h_B(b)$,
but we are sure that the corresponding $L$-interpretations $(h_A \circ \pi_A)$ and $(h_B \circ \pi_B)$
are elementarily equivalent, then we can exclude $(a, b)$ as a witness for $\pi_A \not\equiv \pi_B$.
If we are able to provide enough pairs of homomorphisms so that each distinct pair $(a, b)$ can be excluded,
then $\pi_A \equiv \pi_B$ must hold. The following definition formalises the required properties.

\begin{Definition}[Separating Homomorphism Pairs]
A set $S \subseteq \Hom^2(K, L)$ of homomorphism pairs $h_A, h_B \from K \to L$ is called \emph{separating}
if for all $a, b \in K$ with $a \neq b$, there is a pair $(h_A, h_B) \in S$ such that $h_A(a) \neq h_B(b)$.
$S$ is called \emph{diagonal} if $h_A = h_B$ for all pairs $(h_A, h_B) \in S$.
In that case, we may write $S$ as a subset of $\Hom(K, L)$.
\end{Definition}

Note that a single \emph{injective} homomorphism $h \from K \to L$ induces
the diagonal separating set $S \coloneqq \{(h, h)\}$ with just one element.
Moreover, some semirings, such as $\PosBool[X]$, can be completely decomposed into $L \coloneqq \Bool$
using a diagonal separating set of semiring homomorphisms as follows.
Any subset $Y\subseteq X$ induces a unique homomorphism $h_Y \from \PosBool[X] \to \Bool$ by
$h_Y(x) = \top$ for $x \in Y$ and $h_Y(x) = \bot$ for $x\in X\setminus Y$.
Clearly, for any $p \in \PosBool[X]$, we have that $h_Y(p) = \top$ if, and only if,
$p$ contains a monomial with only variables from $Y$.

\begin{lemma}
\label{lem:separating}
The set $S:=\{h_Y \mid Y \subseteq X\} \subseteq \Hom(\PosBool[X], \Bool)$ is a diagonal separating set of homomorphisms.
\end{lemma}

\begin{proof}
Consider $p, q \in \PosBool[X]$ such that $p\neq q$.
Among the monomials that appear in one of the two polynomials $p,q$
but not in the other, let $m$ be one whose set $Y$ of variables is minimal.
By symmetry, we can assume that $m$ appears in $p$ but not in $q$.
It follows that $h_Y(p)=\top$. We claim that $h_Y(q)=\bot$.
Otherwise $q$ must contain a monomial $m'$ with only variables from $Y$.
Since $m'$ has less variables than $m$, $m'$ must also be contained in $p$.
But $m'$ absorbs $m$, so $m$ does not occur in $p$, a contradiction.
\end{proof}

On the other side, $\Nat[X]$, among other semirings, cannot be decomposed into $\Bool$ by a diagonal separating set. 
For example, $h(x + xy) = h(x) \lor (h(x) \land h(y)) = h(x)$ for all homomorphisms
$h \from \Nat[X] \to \Bool$, but $x + xy \neq x$.
The reason why a decomposition into $\Bool$ would be useful for model theory is given by the
following reduction technique.

\begin{proposition}[Reduction Technique]
Let $\pi_A \from \Lit_A(\tau) \to K$ and $\pi_B \from \Lit_B(\tau) \to K$ be two $K$-interpretations,
$\ta \in A^k$ and $\tb \in B^k$ be $k$-tuples and $S \subseteq \Hom^2(K, L)$ a separating set of homomorphism pairs.
Then, for any formula $\phi(x_1, \dots, x_k) \in \FO(\tau)$, we have that whenever
$(h_A \circ \pi_A) \intpr{\phi(\ta)} = (h_B \circ \pi_B) \intpr{\phi(\tb)}$ for all $(h_A, h_B) \in S$,
then also $\pi_A \intpr{\phi(\ta)} = \pi_B \intpr{\phi(\tb)}$.
\end{proposition}

\begin{proof}
Suppose that $\pi_A \intpr{\phi(\ta)} \neq \pi_B \intpr{\phi(\tb)}$.
Then, by definition of $S$, there exists a pair $(h_A, h_B) \in S$ such that
$h_A(\pi_A \intpr{\phi(\ta)}) \neq h_B(\pi_B \intpr{\phi(\tb)})$.
Applying \hyperref[lem:fundamental]{the fundamental property} yields
$(h_A \circ \pi_A) \intpr{\phi(\ta)} \neq (h_B \circ \pi_B) \intpr{\phi(\tb)}$.
\end{proof}

\begin{corollary}
\label{cor:reduction}
With $S$ as above, $(h_A \circ \pi_A), \ta \equiv (h_B \circ \pi_B), \tb$ for all $(h_A, h_B) \in S$
implies that $\pi_A, \ta \equiv \pi_B, \tb$. Moreover, for each $m \in \Nat$,
$(h_A \circ \pi_A), \ta \equiv_m (h_B \circ \pi_B), \tb$ for all $(h_A, h_B) \in S$
implies $\pi_A, \ta \equiv_m \pi_B, \tb$.
\end{corollary}

If we choose the target semiring $L \coloneqq \Bool$, then \hyperref[cor:reduction]{the corollary} shows that
proving equivalence in $K$ may be reduced to proving equivalence in $\Bool$, which allows us to re-use results from standard semantics.

\subsection{Applications to $\PosBool[X]$ and $\Why[X]$}

Consider the following two $\PosBool[X]$-interpretations $\pi_{xy}, \pi_{yx}$ with $X \coloneqq \{x, y\}$
over the universe $A \coloneqq \{a, b, c, d\}$ with four elements and a signature
$\tau \coloneqq \{P, Q\}$ with two unary relation symbols.

\vspace{\baselineskip}
\begin{minipage}{\linewidth}
\centering
$\pi_{xy}:\quad$
\begin{tabular}{c | c | c | c | c |}
$A$ & $P$ & $Q$ & $\neg P$ & $\neg Q$ \\ \hline
$a$ & $0$ & $y$ & $x$ & $0$ \\
$b$ & $x$ & $0$ & $0$ & $y$ \\
$c$ & $y$ & $x$ & $0$ & $0$ \\
$d$ & $0$ & $0$ & $y$ & $x$ \\
\end{tabular}
$\quad\quad\quad\pi_{yx}:\quad$
\begin{tabular}{c | c | c | c | c |}
$A$ & $P$ & $Q$ & $\neg P$ & $\neg Q$ \\ \hline
$a$ & $y$ & $0$ & $0$ & $x$ \\
$b$ & $0$ & $x$ & $y$ & $0$ \\
$c$ & $x$ & $y$ & $0$ & $0$ \\
$d$ & $0$ & $0$ & $x$ & $y$ \\
\end{tabular}\hspace*{8mm}
\end{minipage}
\vspace{\baselineskip}

Thanks to \hyperref[lem:separating]{Lemma~\ref*{lem:separating}},
the four homomorphisms $S = \{h_{\emptyset}, h_{\{x\}}, h_{\{y\}}, h_{\{x, y\}}\}$
induce a separating set on $\PosBool[X]$ and with \hyperref[cor:reduction]{Corollary~\ref*{cor:reduction}},
it suffices to show that $(h \circ \pi_{xy}) \equiv (h \circ \pi_{yx})$ in $\Bool$ for all $h \in S$
in order to prove $\pi_{xy} \equiv \pi_{yx}$ on $\PosBool[X]$.
Indeed, the following tables demonstrate that $(h \circ \pi_{xy}) \iso (h \circ \pi_{yx})$ holds for all $h \in S$.

\vspace{\baselineskip}
\begin{minipage}{\linewidth}
\centering
$\phantom{h_{\{y\}} \circ \pi_{xy}:}\mathllap{h_{\emptyset} \circ \pi_{xy}:}\;$
\begin{tabular}{c | c | c | c | c |}
$A$ & \ooalign{\hphantom{\colorbox{black!15}{$\top$}}\cr\hfil$P$\hfil\cr} & \ooalign{\hphantom{\colorbox{black!15}{$\top$}}\cr\hfil$Q$\hfil\cr} & $\neg P$ & $\neg Q$ \\ \hline
$a$ & $\bot$ & $\bot$ & $\bot$ & $\bot$ \\
$b$ & $\bot$ & $\bot$ & $\bot$ & $\bot$ \\
$c$ & $\bot$ & $\bot$ & $\bot$ & $\bot$ \\
$d$ & $\bot$ & $\bot$ & $\bot$ & $\bot$ \\
\end{tabular}
$\quad\phantom{h_{\{y\}} \circ \pi_{yx}:}\mathllap{h_{\emptyset} \circ \pi_{yx}:}\;$
\begin{tabular}{c | c | c | c | c |}
$A$ & \ooalign{\hphantom{\colorbox{black!15}{$\top$}}\cr\hfil$P$\hfil\cr} & \ooalign{\hphantom{\colorbox{black!15}{$\top$}}\cr\hfil$Q$\hfil\cr} & $\neg P$ & $\neg Q$ \\ \hline
$a$ & $\bot$ & $\bot$ & $\bot$ & $\bot$ \\
$b$ & $\bot$ & $\bot$ & $\bot$ & $\bot$ \\
$c$ & $\bot$ & $\bot$ & $\bot$ & $\bot$ \\
$d$ & $\bot$ & $\bot$ & $\bot$ & $\bot$ \\
\end{tabular}\hspace*{8mm}
\end{minipage}

\vspace{\baselineskip}

\begin{minipage}{\linewidth}
\centering
$\phantom{h_{\{y\}} \circ \pi_{xy}:}\mathllap{h_{\{x\}} \circ \pi_{xy}:}\;$
\begin{tabular}{c | c | c | c | c |}
$A$ & \ooalign{\hphantom{\colorbox{black!15}{$\top$}}\cr\hfil$P$\hfil\cr} & \ooalign{\hphantom{\colorbox{black!15}{$\top$}}\cr\hfil$Q$\hfil\cr} & $\neg P$ & $\neg Q$ \\ \hline
\vphantom{$X^{X^{X^X}}$}$a$ & $\bot$ & $\bot$ & \colorbox{black!15}{$\top$} & $\bot$ \\
$b$ & \colorbox{black!15}{$\top$} & $\bot$ & $\bot$ & $\bot$ \\
$c$ & $\bot$ & \colorbox{black!15}{$\top$} & $\bot$ & $\bot$ \\
$d$ & $\bot$ & $\bot$ & $\bot$ & \colorbox{black!15}{$\top$} \\
\end{tabular}
$\quad\phantom{h_{\{y\}} \circ \pi_{yx}:}\mathllap{h_{\{x\}} \circ \pi_{yx}:}\;$
\begin{tabular}{c | c | c | c | c |}
$A$ & \ooalign{\hphantom{\colorbox{black!15}{$\top$}}\cr\hfil$P$\hfil\cr} & \ooalign{\hphantom{\colorbox{black!15}{$\top$}}\cr\hfil$Q$\hfil\cr} & $\neg P$ & $\neg Q$ \\ \hline
\vphantom{$X^{X^{X^X}}$}$a$ & $\bot$ & $\bot$ & $\bot$ & \colorbox{black!15}{$\top$} \\
$b$ & $\bot$ & \colorbox{black!15}{$\top$} & $\bot$ & $\bot$ \\
$c$ & \colorbox{black!15}{$\top$} & $\bot$ & $\bot$ & $\bot$ \\
$d$ & $\bot$ & $\bot$ & \colorbox{black!15}{$\top$} & $\bot$ \\
\end{tabular}\hspace*{8mm}
\end{minipage}

\vspace{\baselineskip}

\begin{minipage}{\linewidth}
\centering
$\phantom{h_{\{y\}} \circ \pi_{xy}:}\mathllap{h_{\{y\}} \circ \pi_{xy}:}\;$
\begin{tabular}{c | c | c | c | c |}
$A$ & \ooalign{\hphantom{\colorbox{black!15}{$\top$}}\cr\hfil$P$\hfil\cr} & \ooalign{\hphantom{\colorbox{black!15}{$\top$}}\cr\hfil$Q$\hfil\cr} & $\neg P$ & $\neg Q$ \\ \hline
\vphantom{$X^{X^{X^X}}$}$a$ & $\bot$ & \colorbox{black!15}{$\top$} & $\bot$ & $\bot$ \\
$b$ & $\bot$ & $\bot$ & $\bot$ & \colorbox{black!15}{$\top$} \\
$c$ & \colorbox{black!15}{$\top$} & $\bot$ & $\bot$ & $\bot$ \\
$d$ & $\bot$ & $\bot$ & \colorbox{black!15}{$\top$} & $\bot$ \\
\end{tabular}
$\quad\phantom{h_{\{y\}} \circ \pi_{yx}:}\mathllap{h_{\{y\}} \circ \pi_{yx}:}\;$
\begin{tabular}{c | c | c | c | c |}
$A$ & \ooalign{\hphantom{\colorbox{black!15}{$\top$}}\cr\hfil$P$\hfil\cr} & \ooalign{\hphantom{\colorbox{black!15}{$\top$}}\cr\hfil$Q$\hfil\cr} & $\neg P$ & $\neg Q$ \\ \hline
\vphantom{$X^{X^{X^X}}$}$a$ & \colorbox{black!15}{$\top$} & $\bot$ & $\bot$ & $\bot$ \\
$b$ & $\bot$ & $\bot$ & \colorbox{black!15}{$\top$} & $\bot$ \\
$c$ & $\bot$ & \colorbox{black!15}{$\top$} & $\bot$ & $\bot$ \\
$d$ & $\bot$ & $\bot$ & $\bot$ & \colorbox{black!15}{$\top$} \\
\end{tabular}\hspace*{8mm}
\end{minipage}

\vspace{\baselineskip}

\begin{minipage}{\linewidth}
\centering
$\phantom{h_{\{y\}} \circ \pi_{xy}:}\mathllap{h_{X} \circ \pi_{xy}:}\;$
\begin{tabular}{c | c | c | c | c |}
$A$ & \ooalign{\hphantom{\colorbox{black!15}{$\top$}}\cr\hfil$P$\hfil\cr} & \ooalign{\hphantom{\colorbox{black!15}{$\top$}}\cr\hfil$Q$\hfil\cr} & $\neg P$ & $\neg Q$ \\ \hline
\vphantom{$X^{X^{X^X}}$}$a$ & $\bot$ & \colorbox{black!15}{$\top$} & \colorbox{black!15}{$\top$} & $\bot$ \\
$b$ & \colorbox{black!15}{$\top$} & $\bot$ & $\bot$ & \colorbox{black!15}{$\top$} \\
$c$ & \colorbox{black!15}{$\top$} & \colorbox{black!15}{$\top$} & $\bot$ & $\bot$ \\
$d$ & $\bot$ & $\bot$ & \colorbox{black!15}{$\top$} & \colorbox{black!15}{$\top$} \\
\end{tabular}
$\quad\phantom{h_{\{y\}} \circ \pi_{yx}:}\mathllap{h_{X} \circ \pi_{yx}:}\;$
\begin{tabular}{c | c | c | c | c |}
$A$ & \ooalign{\hphantom{\colorbox{black!15}{$\top$}}\cr\hfil$P$\hfil\cr} & \ooalign{\hphantom{\colorbox{black!15}{$\top$}}\cr\hfil$Q$\hfil\cr} & $\neg P$ & $\neg Q$ \\ \hline
\vphantom{$X^{X^{X^X}}$}$a$ & \colorbox{black!15}{$\top$} & $\bot$ & $\bot$ & \colorbox{black!15}{$\top$} \\
$b$ & $\bot$ & \colorbox{black!15}{$\top$} & \colorbox{black!15}{$\top$} & $\bot$ \\
$c$ & \colorbox{black!15}{$\top$} & \colorbox{black!15}{$\top$} & $\bot$ & $\bot$ \\
$d$ & $\bot$ & $\bot$ & \colorbox{black!15}{$\top$} & \colorbox{black!15}{$\top$} \\
\end{tabular}\hspace*{8mm}
\end{minipage}
\vspace{\baselineskip}

Thus, we conclude that $\pi_{xy} \equiv \pi_{yx}$ and due to $\pi_{xy} \not\iso \pi_{yx}$, this shows that
for finite $\PosBool[X]$-interpretations, elementary equivalence does not necessarily imply isomorphism. 
Moreover, similar examples can be constructed for any distributive lattice semiring $K$ thanks to the universal property of 
$\PosBool[X]$, by assigning $r, s \in K$ to the variables $x, y$:

\vspace{\baselineskip}
\begin{minipage}{\linewidth}
\centering
$\pi_{rs}:\quad$
\begin{tabular}{c | c | c | c | c |}
$A$ & $P$ & $Q$ & $\neg P$ & $\neg Q$ \\ \hline
$a$ & $0$ & $s$ & $r$ & $0$ \\
$b$ & $r$ & $0$ & $0$ & $s$ \\
$c$ & $s$ & $r$ & $0$ & $0$ \\
$d$ & $0$ & $0$ & $s$ & $r$ \\
\end{tabular}
$\quad\quad\quad\pi_{sr}:\quad$
\begin{tabular}{c | c | c | c | c |}
$A$ & $P$ & $Q$ & $\neg P$ & $\neg Q$ \\ \hline
$a$ & $s$ & $0$ & $0$ & $r$ \\
$b$ & $0$ & $r$ & $s$ & $0$ \\
$c$ & $r$ & $s$ & $0$ & $0$ \\
$d$ & $0$ & $0$ & $r$ & $s$ \\
\end{tabular}\hspace*{8mm}
\end{minipage}
\vspace{\baselineskip}

Clearly, $\pi_{rs} \equiv \pi_{sr}$ holds as above, and the only requirement for $\pi_{rs} \not\iso \pi_{sr}$
is that $r$ and $s$ must be distinct. This yields the following theorem.

\begin{theorem}
\label{thm:threeelement}
For any distributive lattice semiring $K$ with at least three elements,
there is a pair of finite $K$-interpretations $\pi_{rs}, \pi_{sr}$
over a universe with four elements and a signature with two unary relation symbols
such that $\pi_{rs} \equiv \pi_{sr}$, but $\pi_{rs} \not\iso \pi_{sr}$.
\end{theorem}

Note that the two $K_4$-interpretations $\pi_{PQ}$ and $\pi_{QP}$
from the opening example of \hyperref[sec:separation]{this section}
can be shown to be elementarily equivalent using a similar technique as above.
In fact, \hyperref[thm:threeelement]{the above theorem} even shows that the opening example was not minimal
and a counterexample with only three semiring elements in $K_3 = \{0, 1, 2\}$ exists.

We shall strengthen the result from \hyperref[thm:threeelement]{Theorem~\ref*{thm:threeelement}}
to the class of all fully idempotent semirings by simply regarding $\pi_{xy}$ and $\pi_{yx}$
as $\Why[X]$-interpretations instead of $\PosBool[X]$-interpretations.
However, the proof that $\pi_{xy} \equiv \pi_{yx}$ becomes more involved,
since a diagonal separating set for $\Why[X]$ into $\Bool$ does not exist.
Nevertheless, a suitable separating set can be obtained by exploiting homomorphisms into $\Why[X]$ itself.
Consider any permutation $\sigma \from X \to X$ of the variables.
Surely, it induces an automorphism $h_{\sigma}$ of $\Why[X]$. 
In the previous example, if $\sigma$ swaps the variables $x$ and $y$,
then applying $h_{\sigma}$ to $\pi_{xy}$ yields an interpretation
that is isomorphic to $\pi_{yx}$, as illustrated below.

\vspace{\baselineskip}
\begin{minipage}{\linewidth}
\centering
$(h_{\sigma} \circ \pi_{xy}):\quad$
\begin{tabular}{c | c | c | c | c |}
$A$ & $P$ & $Q$ & $\neg P$ & $\neg Q$ \\ \hline
$a$ & $0$ & $x$ & $y$ & $0$ \\
$b$ & $y$ & $0$ & $0$ & $x$ \\
$c$ & $x$ & $y$ & $0$ & $0$ \\
$d$ & $0$ & $0$ & $x$ & $y$ \\
\end{tabular}
$\quad\quad\quad\pi_{yx}:\quad$
\begin{tabular}{c | c | c | c | c |}
$A$ & $P$ & $Q$ & $\neg P$ & $\neg Q$ \\ \hline
$a$ & $y$ & $0$ & $0$ & $x$ \\
$b$ & $0$ & $x$ & $y$ & $0$ \\
$c$ & $x$ & $y$ & $0$ & $0$ \\
$d$ & $0$ & $0$ & $x$ & $y$ \\
\end{tabular}\hspace*{8mm}
\end{minipage}
\vspace{\baselineskip}

With this insight, we can construct a suitable separating set $S \subseteq \Hom^2(\Why[X], \Why[X])$,
starting with the pair $(h_{\sigma}, h_{\id}) \in S$. This pair alone does not separate $\Why[X]$,
since we have $x \neq y$ in $\Why[X]$, but $h_{\sigma}(x) = y = h_{\id}(y)$.
Hence, we add more homomorphisms by annihilating some variables,
similarly to the construction for $\PosBool[X]$ in \hyperref[lem:separating]{Lemma~\ref*{lem:separating}}.
Fixing a permutation $\sigma \from X \to X$, we want to construct a homomorphism $h_{\sigma}^Y$
that annihilates all variables in $X \setminus Y$ and permutes the variables in $Y$.
Observe that for each $x\in Y$ there is a minimal number $r(x)\geq 1$ such that $\sigma^{r(x)}(x)\in Y$.
Formally, we define $\sigma^Y \from X \to Y \cup \{0\}$ by setting $\sigma^Y(x) \coloneqq \sigma^{r(x)}(x)$ for $x\in Y$,
and $\sigma^Y(x) \coloneqq 0$ for $x\in X\setminus Y$.
Note that $\sigma^Y$ induces a homomorphism $h_{\sigma}^Y \from \Why[X] \to \Why[X]$.

\begin{lemma}
$S \coloneqq \{(h_{\sigma}^Y, h_{\id}^Y) \mid Y \subseteq X\} \subseteq \Hom^2(\Why[X], \Why[X])$ is a separating set of homomorphism pairs.
\end{lemma}

\begin{proof}
Suppose that $p \neq q$ for a pair $p, q \in \Why[X]$.
A monomial in $\Why[X]$ can be identified with the set of its variables.
Thus, without loss of generality, there is some $Y\subseteq X$ with $Y \in p$ and $Y \notin q$. 
Surely, $h_{\sigma}^Y(p)$ contains the monomial $h_{\sigma}^Y(Y) = Y$,
but $h_{\id}^Y(q)$ only contains monomials from $q$,
hence it does not contain $Y$ and $h_{\sigma}^Y(p) \neq h_{\id}^Y(q)$.
\end{proof}

Before applying the reduction technique to obtain $\pi_{xy} \equiv \pi_{yx}$
from \hyperref[cor:reduction]{Corollary~\ref*{cor:reduction}}, it only remains to show
that $(h_{\sigma}^Y \circ \pi_{xy}) \equiv (h_{\id}^Y \circ \pi_{yx})$ for all $Y \subseteq X \coloneqq \{x, y\}$.
Since we have already illustrated $(h_{\sigma}^X \circ \pi_{xy}) \iso (h_{\id}^X \circ \pi_{yx})$ above,
we only need to consider the cases where $Y \subsetneq X$.
But then, at most one variable is contained in $Y$ and the remaining variables are annihilated
by $h_{\sigma}^Y$ and $h_{\id}^Y$, thus $(h_{\sigma}^Y \circ \pi_{xy}) \iso (h_{\id}^Y \circ \pi_{yx})$ clearly follows.
The reduction technique then implies that $\pi_{xy} \equiv \pi_{yx}$ on $\Why[X]$,
which can naturally be lifted to all fully idempotent semirings thanks to \hyperref[lem:universal]{the universal property}.

\begin{theorem}
For any fully idempotent semiring $K$ with at least three elements,
there is a pair of finite $K$-interpretations $\pi_{rs}, \pi_{sr}$ over a universe with four elements
and a signature with two unary relation symbols such that
$\pi_{rs} \equiv \pi_{sr}$, but $\pi_{rs} \not\iso \pi_{sr}$.
\end{theorem}

In conclusion, the proof of $\pi_{xy} \equiv \pi_{yx}$ on $\PosBool[X]$ illustrates
how elementary equivalence in semiring semantics can be reduced
to elementary equivalence on $\Bool$ by completely decomposing
$\PosBool[X]$ to $\Bool$ with homomorphisms. Moreover,
the proof of elementary equivalence on $\Why[X]$ shows that
it may even pay off to use separating sets of homomorphisms
from $\Why[X]$ to $\Why[X]$ itself.

\section{Characteristic sentences}

The results of the previous section raise the question whether it is possible to construct a similar example
of non-isomorphic, but elementarily equivalent interpretations also for the most general semiring $\Nat[X]$,
and lift it to all commutative semirings with at least three elements.
In order to show that this is not possible, we draw inspiration from classical semantics,
where for each finite $\tau$-structure $\StrA$ with universe $A=\{a_1, \dots, a_n\}$ 
one can construct a \emph{characteristic sentence} $\chi_{\StrA}$
such that $\StrB \models \chi_{\StrA}$ if, and only if, $\StrA \iso \StrB$. 
The characteristic sentence is explicitly defined as
\begin{alignat*}{1}
\chi_{\StrA} &\coloneqq \exists x_1 \ldots \exists x_n (\phi(\tx) \land \psi(\tx)) \quad \text{with}\\
\phi(\tx) &\coloneqq \bigwedge_{1 \le i < j \le n} x_i \neq x_j \land \forall y \bigvee_{i \le n} y = x_i \quad\text{and} \quad
\psi(\tx) \coloneqq \bigwedge \{L(\tx) \in \Lit_n(\tau) \mid \StrA \models L(\ta)\}.
\end{alignat*}
The subformula $\phi(\tx)$ of this sentence asserts that the universe has precisely $n$ elements 
assigned to the variables $\tx$. Since $\phi(\tx)$ uses only equalities and inequalities 
it can be used as-is for semiring semantics in any semiring.

\begin{lemma}
\label{lem:samesize}
For every $K$-interpretation $\pi_B \from \Lit_B(\tau) \to K$ into an arbitrary semiring $K$
and every tuple $\tb = (b_1,\dots,b_n)$, we have that $\pi_B \intpr{\phi(\tb)} = 1$
if $B = \{b_1, \dots, b_n\}$ and $b_i\neq b_j$ for $i\neq j$,
and $\pi_B \intpr{\phi(\tb)} = 0$, otherwise.
\end{lemma}

\begin{proof}
Semiring interpretations evaluate equalities and inequalities to 0 and 1, so
\[\pi_B \intpr{\phi(\tb)} = \prod_{i<j} \pi_B \intpr{b_i \neq b_j} \cdot \prod_{b \in B}\bigl(\sum_{i\leq n} \pi_B \intpr{b = b_i}\bigr)\]
evaluates to 1 if $b_1, \dots, b_n$ is a distinct enumeration of all elements of $B$, and to 0 otherwise.
\end{proof}

On the other side, $\psi(\ta)$ is the conjunction of all true literals in $\StrA$.
Since $\StrA$ satisfies precisely one literal out of each pair of opposing literals $L$ and $\Lcomp$,
it is clear that $\psi(\ta)$ describes $\StrA$ up to isomorphism.
However, this approach does not lift to arbitrary semiring interpretations $\pi_A$, since different literals in $\pi_A$ 
may have different non-zero values, but conjunctions are interpreted as products, and it is in general 
impossible to trace the result back to the contributions of the literals.

\subsection{The Viterbi semiring}
 
The Viterbi semiring $\Vit = ([0, 1]_{\Real}, \max, \cdot, 0, 1)$ is used in confidence analysis, 
probabilistic parsing, and Hidden Markov Models
(see \cite{DrosteKui09,Goodman99}). It is isomorphic to the tropical semiring $\Trop= (\mathbb{R}_{+}^{\infty},\min,+,\infty,0)$, used for instance for cost analysis and performance evaluation, via $x\mapsto e^{-x}$. 
Hence, all results that we establish for the Viterbi semiring also hold for the tropical semiring.
We can illustrate the shortcomings of the characteristic sentences in their classical form by very simple $\Vit$-interpretations with one element.

\vspace{\baselineskip}
\begin{minipage}{\linewidth}
\centering
$\pi_{19}:\quad$
\begin{tabular}{c | c | c | c | c |}
$A$ & $P$ & $Q$ & $\neg P$ & $\neg Q$ \\ \hline
$a$ & $0.1$ & $0.9$ & $0$ & $0$ \\
\end{tabular}
$\quad\quad\quad\pi_{91}:\quad$
\begin{tabular}{c | c | c | c | c |}
$A$ & $P$ & $Q$ & $\neg P$ & $\neg Q$ \\ \hline
$a$ & $0.9$ & $0.1$ & $0$ & $0$ \\
\end{tabular}
\end{minipage}
\vspace{\baselineskip}

They are clearly not isomorphic, but trying to construct $\chi_{19}$ from $\pi_{19}$ as above
would yield $\chi_{19} = \exists x (\phi(x) \land \psi(x))$ with $\psi(x) = Px \land Qx$,
hence $\pi_{19} \intpr{\chi_{19}} = 0.1 \cdot 0.9 = 0.9 \cdot 0.1 = \pi_{91} \intpr{\chi_{19}}$. 

\medskip

However, under semiring semantics, and especially on the Viterbi semiring $\Vit$,
multiplication need not be idempotent; hence we can hope to distinguish two interpretations
by simply repeating one of the literals.
In the given example, we can set $\psi(x) \coloneqq Px \land (Qx)^2$,
which is short for $Px \land Qx \land Qx$, to obtain
$\pi_{19} \intpr{\chi_{19}} = 0.1 \cdot 0.9^2 \neq 0.9 \cdot 0.1^2 = \pi_{91} \intpr{\chi_{19}}$. 
We now generalise this idea to arbitrary finite $\Vit$-interpretations.
We shall associate with every finite $\Vit$-interpretation $\pi_A \from \Lit_A(\tau) \to \Vit$ and every
$\epsilon \in \Real^+$ a characteristic sentence 
\[\chi_{\pi_A, \epsilon} \coloneqq \exists x_1 \ldots \exists x_n (\phi(\tx) \land \psi_{\epsilon}(\tx)),\] 
with $n \coloneqq \vert A \vert$ and $\phi(x)$ as before, but
a more involved construction of $ \psi_{\epsilon}(\tx)$:
 
Let $\ta = (a_1, \ldots, a_n)$ be some fixed order on $A$ and $L_1(\ta), \ldots, L_k(\ta)$
an arbitrary enumeration of the ``true'' literals in $\Lit_A(\tau)$ with $\pi_A(L_i(\ta)) \neq 0$.
Further, fix a sequence $f(1), \dots, f(k)$ of ``exponents'' in $\N$, where $f(1)=1$ and
$f(i+1)$ is chosen large enough so that 
\[(*) \quad\quad\quad (1-\epsilon)^{f(i+1)} < \epsilon^{f(1)+\dots + f(i)}.\]
Then, put $\psi_{\epsilon}(\tx) \coloneqq \bigwedge_{i = 1}^k L_i(\tx)^{f(i)}$,
where ``exponentiation'' denotes repetition of a literal.

The idea is that $\chi_{\pi_A, \varepsilon}$ should characterise $\pi_A$ up to isomorphism
by repeating the literals in $\psi_{\varepsilon}(\tx)$ ``sufficiently often''
so that the contribution of each literal can be distinguished and
changing the value for one literal surely alters the final value of $\psi_{\varepsilon}(\tx)$.
Since the elements of $\Vit$ are from $[0, 1]_{\Real}$, the values of the literals
can change by an arbitrarily small amount, hence the ``exponents'' $f(i)$ must 
depend on the ``smallest possible change'' $\varepsilon$. This intuition is formalised as follows.

\begin{proposition}
\label{prop:epscharacteristic}
Let $\pi_A \from \Lit_A(\tau) \to \Vit$ and $\pi_B \from \Lit_B(\tau) \to \Vit$
be two finite, model-defining $\Vit$-interpretations, which induce the finite set of values
\[V \coloneqq \{\pi_A(L) \mid L \in \Lit_A(\tau)\} \cup \{\pi_B(L) \mid L \in \Lit_B(\tau)\}.\]
Then, for every $\epsilon \in \Real$ bounded by $0 < \epsilon \le \min \{\vert r - s \vert \mid r, s \in V, r \neq s\}$,
we have that $\pi_A \intpr{\chi_{\pi_A, \epsilon}} = \pi_B \intpr{\chi_{\pi_A, \epsilon}}$ implies $\pi_A \iso \pi_B$.
\end{proposition}

\begin{proof}
Assume $\pi_A \intpr{\chi_{\pi_A, \epsilon}} = \pi_B \intpr{\chi_{\pi_A, \epsilon}}$. 
By construction, $\pi_A \intpr{\chi_{\pi_A, \epsilon}} > 0$, so $\pi_B \intpr{\chi_{\pi_A, \epsilon}} > 0$ as well.
By \hyperref[lem:samesize]{Lemma~\ref*{lem:samesize}}, together with the fact that the existential quantifiers $\exists x_1 \ldots \exists x_n$
in $\chi_{\pi_A, \epsilon}$ are interpreted as $\max$ in the Viterbi semiring $\Vit$, 
this implies that $\vert A \vert = \vert B \vert$, and that we have enumerations
$\ta = (a_1, \ldots, a_n)$ and $\tb = (b_1, \ldots, b_n)$ of the elements of $A$ and $B$, such that 
\[\pi_A \intpr{\chi_{\pi_A, \epsilon}} = \pi_A \intpr{\psi_{\epsilon}(\ta)} = \pi_B \intpr{\psi_{\epsilon}(\tb)} = \pi_B \intpr{\chi_{\pi_A, \epsilon}}.\]

Recall that $\pi_A \intpr{\psi_{\epsilon}(\ta)} = \prod_{i = 1}^k \pi_A(L_i(\ta))^{f(i)} > 0$,
hence $\pi_A(L_1(\ta)), \dots, \pi_A(L_k(\ta))$ are all positive.
Accordingly, $\pi_B \intpr{\psi_{\epsilon}(\tb)} = \prod_{i = 1}^k \pi_B(L_i(\tb))^{f(i)} > 0$, so that
$\pi_B(L_1(\tb)), \dots, \pi_B(L_k(\tb))$ are positive as well.
Given that $\pi_A$ and $\pi_B$ share the same signature and universe size,
any permutations $\ta$, $\tb$ of their elements yield the same number of positive literals, which is $k$ by definition.
We infer that all remaining literals in both interpretations are mapped to zero.
Hence, for $i=1,\dots, k$, let $r_i:= \pi_A(L_i(\ta)) > 0$ and $s_i:= \pi_B(L_i(\tb)) > 0$ be the values of the positive literals,
then it only remains to show that $r_i=s_i$ for all $ i \le k$ in order to conclude 
that $\ta \mapsto \tb$ is indeed an isomorphism from $\pi_A$ to $\pi_B$.

Towards a contradiction, assume that this is not the case
and let $j$ be the maximal index among $1, \dots, k$ 
with $r_j\neq s_j$. 
We can assume that $r_j < s_j$. Since the difference
between the two values is at least $\epsilon$ and since $s_j\leq 1$
it follows that $r_j\leq s_j-\epsilon\leq s_j -\epsilon s_j=(1-\epsilon)s_j$.
Further, we have $\epsilon\leq s_i, r_i\leq 1$ for all $i$.
It follows that
\[ r_1^{f(1)}\cdots r_j^{f(j)}\leq r_j^{f(j)}\leq (1-\epsilon)^{f(j)}s_j^{f(j)} \; \stackrel{*}{<} \; \epsilon^{f(1)+\dots+f(j-1)}\cdot s_j^{f(j)}
\leq s_1^{f(1)}\cdots s_j^{f(j)}.\]
However, since $r_i = s_i$ for $i = j+1, \dots, k$, this would imply that
\[\pi_A \intpr{\psi_{\epsilon}(\ta)} = \prod_{i \leq k} r_i^{f(i)} \neq \prod_{i \leq k} s_i^{f(i)} = \pi_B \intpr{\psi_{\epsilon}(\tb)}\]
and hence $\pi_A \intpr{\chi_{\pi_A, \epsilon}} \neq \pi_B \intpr{\chi_{\pi_A, \epsilon}}$.
\end{proof}

Notice that none of the sentences $\chi_{\pi_A, \epsilon}$ characterises $\pi_A$ alone,
but the countable set $X_{\pi_A} \coloneqq \{\chi_{\pi_A, \epsilon} \mid \epsilon \in \mathbb{Q}^+\}$ does so.
No infinite $\Vit$-interpretation $\pi_B$ agrees with $\pi_A$ on any of the $\epsilon$-characteristic sentences
$\chi_{\pi_A, \epsilon}$ due to $\phi(\tx)$, whereas for each finite $\Vit$-interpretation $\pi_B$,
one can calculate an $\epsilon \in \mathbb{Q}$ to apply the proposition just proved.

\begin{theorem}
For finite $\Vit$-interpretations $\pi_A$ and $\pi_B$, $\pi_A \equiv \pi_B$ implies $\pi_A \iso \pi_B$.
\end{theorem}

As a consequence, there are indeed interesting semirings beyond the Boolean semiring $\Bool$, 
where elementary equivalence implies isomorphism on finite interpretations.

\subsection{Finite axiomatisability}

The characteristic set $X_{\pi_A}$ raises the question whether a finite set of sentences suffices to characterise 
a $\Vit$-interpretation $\pi_A$. We will answer this question positively using two observations.
By \hyperref[prop:epscharacteristic]{Proposition~\ref*{prop:epscharacteristic}}, we observe 
that $\chi_{\pi_A, \epsilon}$ characterises $\pi_A$ up to isomorphism inside the
class of $\Vit$-interpretations that only use values in $V = \{\pi_A(L) \mid L \in \Lit_A(\tau)\}$,
with $\epsilon \coloneqq \min\{\vert r - s\vert \mid r, s \in V, r \neq s\}$.
Hence, $\pi_A$ can be characterised by adding sentences to ensure that no values outside of $V$ are used.

We will show that this is possible by building sentences that fix particular values $\pi_A(L)$.

\begin{Definition}
Let $\pi \from \Lit_A(\tau)\to \Vit$ be a finite $\Vit$-interpretation over a universe $A$ with $n$ elements
and $\phi(\tx) \in \FO(\tau)$ a formula with $k \in \Nat$ free variables.
The sequence $(s_{\pi, \phi}^i)_{1 \le i \le n^k}$ is defined as
the non-increasingly sorted sequence of the values $\pi \intpr{\phi(\ta)}$ for $\ta \in A^k$.
\end{Definition}

In particular, $s_{\pi, \phi}^1$ is the largest possible value $\pi \intpr{\phi(\ta)}$; further,
$s_{\pi, \phi}^2 \le s_{\pi, \phi}^1$ is either the second largest one, or equal to $s_{\pi, \phi}^1$
if the maximal value is shared by two distinct tuples $\ta, \tb \in A^k$, and so on. 
We construct a series of sentences that fix the values $(s_{\pi, \phi}^i)_{1 \le i \le n^k}$.

\begin{lemma}[Sorting Lemma]
\phantomsection\label{lem:sorting}
For $\phi(\tx)\in\FO(\tau)$ with $k$ free variables and $1 \le i \le n^k$, let
\[\psi_{\phi}^i \coloneqq \exists \tx_1 \ldots \exists \tx_i \Bigl(\bigwedge_{1 \le j < \ell \le i} \tx_j \neq \tx_{\ell} \land \bigwedge_{j = 1}^i \phi(\tx_j)\Bigr)\]
where $\tx_1, \ldots \tx_i$ are $k$-tuples of variables.
Then, for any $\Vit$-interpretation $\pi \from \Lit_A(\tau)\to \Vit$
over a universe with $n$ elements, we have that
\[\pi \intpr{\psi_{\phi}^i} = \prod_{1 \le j \le i} s_{\pi, \phi}^j \quad\text{for } 1 \le i \le n.\]
\end{lemma}

\begin{proof}
Recall that existential quantifiers are interpreted as $\max$ on $\Vit$.
Due to monotonicity of multiplication, the maximum $\pi \intpr{\psi_{\pi, \phi}^i}$ is achieved
by picking the $i$ pairwise distinct tuples $\ta_1, \ldots, \ta_i$ that yield the largest values
$\pi \intpr{\phi(\ta_j)} = s_{\pi, \phi}^j$ and inserting them for $\tx_1, \ldots, \tx_i$.
Clearly, this yields $\pi \intpr{\psi_{\phi}^i} = \prod_{1 \le j \le i} s_{\pi, \phi}^j$.
\end{proof}

By observing that $\Vit$ is cancellative, i.e. $ab = ac$ implies $b = c$ for all $a, b, c \in \Vit$ with $a \neq 0$,
we may disentangle the products $\prod_{1 \le j \le i} s_{\pi, \phi}^j$ to draw the following conclusion.

\begin{corollary}
\label{cor:sorting}
Let $\phi(\tx)$ be as above, and consider two $\Vit$-interpretations $\pi \from \Lit_A(\tau)\to \Vit$ and
$\pi' \from \Lit_B(\tau) \to \Vit$ with $\vert A \vert = \vert B \vert = n$. If $\pi$ and $\pi'$ agree on
$\Psi \coloneqq \{\psi_{\phi}^i \mid 1 \le i \le n^k\}$,
then $s_{\pi, \phi}^i = s_{\pi', \phi}^i$ for all $1 \le i \le n^k$.
In other words, the values $\pi \intpr{\phi(\ta)}$ for $\ta \in A^k$ and
$\pi' \intpr{\phi(\tb)}$ for $\tb \in B^k$ are the same, up to permutation.
\end{corollary}

If $R \in \tau$ is a $k$-ary relation, pick $\phi(\tx) \coloneqq R\tx$ and construct
$\Psi_R \coloneqq \{\psi_{R\tx}^i \mid 1 \le i \le n^k\}$ according to \hyperref[lem:sorting]{the Sorting Lemma}.
Similarly, construct $\Psi_{\neg R}$ from $\phi(\tx) \coloneqq \neg R\tx$. Then, define
\[\Psi_{\tau} \coloneqq \bigcup_{R \in \tau} \Psi_R \cup \bigcup_{R \in \tau} \Psi_{\neg R}.\]
Clearly, any two $\Vit$-interpretations over $\tau$ that agree on $\Psi_{\tau}$
use the same set of values from $\Vit$. Putting this together
with the characteristic sentences $\chi_{\pi_A, \epsilon}$
from \hyperref[prop:epscharacteristic]{Proposition~\ref*{prop:epscharacteristic}}
provides a finite axiomatisation of any $\Vit$-interpretation.

\begin{theorem}
Let $\pi \from \Lit_A(\tau) \to \Vit$ be a finite $\Vit$-interpretation.
Then, $\Psi_{\tau} \cup \{\chi_{\pi, \epsilon}\}$ is a finite axiomatisation of $\pi$ up to isomorphism,
where $\epsilon \coloneqq \min \{\vert r - s \vert \mid r, s \in \pi(\Lit_A(\tau)), r \neq s\}$.
\end{theorem}

\begin{proof}
Let $\pi' \from \Lit_B(\tau) \to \Vit$ agree with $\pi$ on all sentences in $\Psi_{\tau} \cup \{\chi_{\pi, \epsilon}\}$.
Due to the construction of $\chi_{\pi, \epsilon}$, $\pi'$ is finite and $\vert A \vert = \vert B \vert$.
Since $\pi'$ agrees with $\pi$ on $\Psi_{\tau}$, we have $\{\pi(L) \mid L \in \Lit_A(\tau)\} = \{\pi'(L) \mid L \in \Lit_B(\tau)\}$.
Thus, we can invoke \hyperref[prop:epscharacteristic]{Proposition~\ref*{prop:epscharacteristic}}
by observing that $V = \pi(\Lit_A(\tau))$ and conclude that
$\pi \intpr{\chi_{\pi, \epsilon}} = \pi' \intpr{\chi_{\pi, \epsilon}}$ implies $\pi \iso \pi'$.
\end{proof}

Under classical semantics, any finite axiom system $\Phi \subseteq \FO(\tau)$ can be collapsed 
to a single axiom $\psi \coloneqq \bigwedge \Phi$, but this is not the case in semiring semantics.
To illustrate this, we shall show that there are $\Vit$-interpretations that cannot be
axiomatised up to isomorphism by a single sentence.

\begin{proposition}
There exist $\Vit$-interpretations $\pi \from \Lit_A(\tau) \to \Vit$ such that, for every sentence $\psi\in\FO(\tau)$,
there exists an interpretation $\pi' \from \Lit_A(\tau) \to \Vit$ such that $\pi \not \iso \pi'$, but
$\pi \intpr{\psi} = \pi' \intpr{\psi}$.
\end{proposition}

\begin{proof}
Take an interpretation with just two atoms $Pa$ and $Qa$ and with values $\pi(Pa)=p$
and $\pi(Qa) = q$ such that $0< p,q <1$ are multiplicatively independent real numbers, i.e. there are no solutions
$k, \ell \in{\mathbb Z}$ to the equation $p^kq^\ell = 1$, except $k = \ell = 0$. 
Let $\pi_{\Bool}$ be the corresponding $\Bool[x,y]$-interpretation, with $\pi_{\Bool}(Pa) = x$
and $\pi_{\Bool}(Qa) = y$.
A sentence $\psi\in\FO$ is evaluated under $\pi_{\Bool}$ to a polynomial 
$\pi_{\Bool}\intpr{\psi}\in\Bool[x,y]$, and by \hyperref[lem:universal]{the universal property} for idempotent semirings,
the homomorphism $h \from \Bool[x,y] \to \Vit$ induced by $h(x) = p$ and $h(y) = q$
maps $\pi_{\Bool} \intpr{\psi}$ to $\pi \intpr{\psi}$.
Writing $\pi_{\Bool}\intpr{\psi}$ as a sum of monomials $m = x^iy^j$, we conclude that
$\pi\intpr{\psi} = p^i q^j$ is the maximal value $m(p,q)$ for the monomials $m$ occuring
in $\pi_{\Bool}\intpr{\psi}$. Since $p,q$ are multiplicatively independent, 
no other monomial can take the same value, i.e. $m'(p,q) < m(p,q)$ for all other monomials 
$m'$ in $\pi_{\Bool}\intpr{\psi}$.
We now can certainly find a value $r\neq p$ that is sufficiently close to $p$,
and a value $s$ such that $r^is^j = p^iq^j$, i.e. $m(r,s) = m(p,q)$, but $m'(r,s) < m(r,s)$
for all other monomials $m'$ in $\pi_{\Bool} \intpr{\psi}$.
For the $\Vit$-interpretation $\pi'$ with $\pi'(Pa) = r$ and $\pi'(Qa) = s$ this implies
that $\pi'\intpr{\psi}=r^i s^j=p^i q^j = \pi \intpr{\psi}$, but clearly, $\pi' \not\iso \pi$.
\end{proof}

This result can be strengthened in many directions. It holds, in fact, for almost all $\Vit$-interpretations,
as long as they do not map all literals to either 0 or 1. Further, we shall exploit the isomorphism
of $\Vit$ and $\Trop$ in order to prove explicit lower bounds
on the number of axioms that are needed to characterise an interpretation,
depending on the number of literals mapped to multiplicatively independent values.

\subsection{Lower bound for axiomatisations of $\Trop$- and $\Vit$-interpretations}

Recall that $\Vit = ([0, 1]_{\Real}, \max, \cdot, 0, 1)$ is isomorphic to $\Trop = (R_{+}^{\infty}, \min, +, \infty, 0)$
via isomorphisms $\vtot(a) = -\log_b(a)$ for any fixed base $b \in \Real_{> 1}$, and
the corresponding inverse isomorphisms $\ttov(a) = b^{-a}$ for any fixed $b \in \Real_{> 1}$.
We formulate our result in terms of $\Trop$.

\begin{theorem}
Let $\pi \from \Lit_A(\tau) \to \Trop$ be any finite, model-defining $\Trop$-interpretation
with $\vert A \vert = n$ and $\vert \Lit_A(\tau) \vert = 2\ell$,
such that its finite values $\pi(\Lit_A(\tau)) \setminus \{\infty\}$ are linearly independent over $\Rat$.
Then, for any set of sentences $\Psi \subseteq \FO(\tau)$ with $\vert \Psi \vert < \ell$,
there is an interpretation $\pi' \from \Lit_A(\tau) \to \Trop$ such that
$\pi \intpr{\psi} = \pi' \intpr{\psi}$ for all $\psi \in \Psi$, but $\pi \not\iso\pi'$.
\end{theorem}

\begin{proof}
Since $\pi$ is model-defining, there are $\ell$ literals $L$ in $\Lit_A(\tau)$ with $\pi(L) \neq \infty$,
which we call the positive literals.
Choose $X \coloneqq \{x_1, \ldots, x_{\ell}\}$ and construct the $\Bool[X]$-interpretation 
$\pi_{\Bool} \from \Lit_A(\tau) \to \Bool[X]$
by assigning a unique variable to each of the positive literals.
Clearly, there is a homomorphism $h \from \Bool[X] \to \Trop$ with $h \circ \pi_{\Bool} = \pi$
induced by mapping each variable to the original value $\pi(L)$ of the corresponding literal.

Enumerate $\Psi = \{\psi_1, \ldots, \psi_{j}\}$ arbitrarily with $j < \ell$ and
construct the polynomials $p_i \coloneqq \pi_{\Bool} \intpr{\psi_i} \in \Bool[X]$ for $1 \le i \le j$.
By \hyperref[lem:fundamental]{the fundamental property},
$\pi \intpr{\psi_i} = h(p_i)$ holds for $1 \le i \le j$.

We assume without loss of generality that $p_i \neq 0$ for all $1 \le i \le j$,
and we will construct $\pi'$ with the same positive literals as $\pi$.
Thus, $p_i$ contains a monomial $m_i$ so that
$h(m_i)$ is minimal among $\{h(m) : m \in p_i\}$.
This monomial is unique thanks to the linear independence of the values of $\pi$,
which guarantees that $h(m) \neq h(m')$ for $m \neq m'$.
Suppose for a contradiction that $h(m) = h(m')$. We may write
\[ h(m) = h \left(\prod_{i = 1}^{\ell} x_i^{m(x_i)}\right)
= \sum_{i = 1}^{\ell} m(x_i) \cdot h(x_i),\]
which implies that $h(m)$ is a linear combination of the values of $\pi$.
In particular, $h(m) = h(m')$ implies that $h(m - m') = 0$,
hence $m(x_i) - m'(x_i) = 0$ for all $1 \le i \le j$ due to linear independence.

We conclude that there is a sufficiently small $\epsilon \in \Real_{>0}$ such that
changing the numbers $h(x_i)$ by less than $\epsilon$ does not affect the monomial order.
In other words, view the values $\overline{v} \coloneqq (h(x_1), \ldots, h(x_{\ell})) \in \Real_{\ge 0}^{\ell}$ as a vector
and notice that any $\overline{w} \in \Real_{\ge 0}^{\ell}$ with $\vert \overline{v} - \overline{w} \vert < \epsilon$
preserves the monomial order, so that if we construct $h' \from \Bool[X] \to \Trop$ induced by
$h'(x_i) = w_i$ for $1 \le i \le \ell$, we have $h(m) < h(m')$ if, and only if, $h'(m) < h'(m')$.

To complete the proof, it remains to ensure that $h'(p_i) = h(p_i)$
stays the same for all $1 \le i \le j$. By the above considerations,
it suffices to ensure that $h'(m_i) = h(m_i)$ for the corresponding
maximal monomials $m_1, \ldots, m_{j}$.
Each of these monomials induces one condition $h(m_i) - h'(m_i) = 0$,
which translates to a linear equation
\begin{alignat*}{1}
h(m_i) - h'(m_i) = \sum_{i = 1}^{\ell} m_i(x_i) (h(x_i) - h'(x_i)) 
= \sum_{i = 1}^{\ell} m_i(x_i) \cdot (v_i - w_i) = 0
\end{alignat*}
on $(\overline{v} - \overline{w})$.

Since there are only $j < \ell$ equations and $\ell$ variables,
the solution space is at least one-dimensional, meaning that we can pick
$\overline{w} \neq \overline{v}$ adequately with $\vert \overline{v} - \overline{w} \vert < \epsilon$
to satisfy all equations and obtain $h'(p_i) = h(p_i)$ for all $1 \le i \le j$.
Note that due to linear independence, none of the entries from $\overline{v}$ was zero,
hence it is possible to ensure that $\overline{w}$ only has positive entries.
We thus can pick $\pi' \coloneqq h' \circ \pi_{\Bool}$ with the desired properties.
\end{proof}

This result translates to $\Vit$ thanks to isomorphism.
Linear independence of values from $\Trop$ as $\Rat$-vectors
translates to multiplicative independence of the corresponding
values from $\Vit$.

\subsection{The semirings $\Nat$ and $\Nat[X]$}

We will now provide a similar analysis of axiomatisablity 
for the most general semiring $\Nat[X]$ by taking a detour via $\Nat$. 
For the construction of the characteristic sentences for $\Nat$, we shall need the following combinatorial lemma.

\begin{lemma}
For any two natural numbers $k, c$ with $c > 1$, there exists a exponent $e$ 
such that, for any two non-decreasing sequences
$r_1 \leq r_2 \leq \dots \leq r_k$ and $s_1 \leq s_2 \leq \dots \leq s_k$ 
of $k$ natural numbers, with $r_k, s_k < c$, the equation
$r_1^e + \cdots + r_k^e=s_1^e + \cdots + s_k^e$ implies that the two
sequences are the same, i.e. $r_i = s_i$ for all $i \leq k$.
\end{lemma}

\begin{proof}
Choose $e$ large enough so that $(c/(c-1))^e > k$. 
Towards a contradiction, assume that there are two \emph{distinct} sequences
$r_1 \leq r_2 \leq \dots r_k < c$ and $s_1 \leq s_2 \leq \dots s_k < c$ 
such that $r_1^e + \cdots + r_k^e = s_1^e + \cdots + s_k^e$.
Let $j$ be the maximal index with $r_j\neq s_j$.
Thanks to additive cancellation, we can remove the summands 
with index $i > j$ to obtain
that $r_1^e + \cdots + r_j^e = s_1^e + \cdots + s_j^e$. 
By symmetry we can assume that $r_j < s_j$.
Since $s_j<c$, it follows that $s_j > (c/(c-1))r_j$ and hence $s_j^e > k\cdot r_j^e$.
But this implies that
\[r_1^e + \cdots + r_j^e \leq j \cdot r_j^e \leq k \cdot r_j^e < s_j^e \leq s_1^e + \cdots + s_j^e,\] 
contradicting the equation above. 
\end{proof}

\begin{lemma}
\label{lem:exponent}
Let $(r_1, \ldots r_k), (s_1, \ldots, s_k) \in \Nat^k$ be strictly bounded by $c$,
that is $r_i, s_i < c$ for all $i \le k$. Then, there is an exponent $e$ depending only on $c$ and $k$ such that
\[\sum_{i = 1}^k r_i^e = \sum_{i = 1}^k s_i^e\]
implies that there is a permutation $\sigma \in S_k$ such that $r_i = s_{\sigma(i)}$ for all $1 \le i \le k$.
\end{lemma}

\begin{proof}
Sort both sequences non-decreasingly, that is, permute them with $\rho, \tau \in S_k$
so that $r_{\rho(1)} \le \ldots \le r_{\rho(k)}$ and $s_{\tau(1)} \le \ldots \le s_{\tau(k)}$. 
By the previous Lemma, there is a suitable $e$ such that $r_{\rho(i)} = s_{\tau(i)}$ for all $i \le k$.
Then, $\sigma \coloneqq \tau \circ \rho^{-1}$ is the desired permutation.
\end{proof}

We are now ready to construct characteristic sentences for 
finite $\Nat$-interpretations $\pi_A \from \Lit_A(\tau) \to \Nat$. For $n=|A|$, let
$L_1(\tx),\dots,L_k(\tx)$ be an enumeration of all literals in $\Lit_n(\tau)$.
For any constant $q\in\Nat$ we define the \emph{$q$-characteristic sentence} $\chi_{\pi_A, q}$ as
\[\chi_{\pi_A, q} \coloneqq \exists x_1 \ldots \exists x_n (\varphi(\tx) \land \psi_q(\tx))^e, \quad\text{with } 
\psi_q(\tx) \coloneqq \bigvee_{i = 1}^k q^{i - 1} \cdot L_i(\tx),\]
with $\varphi(\tx)$ as given before and $e$ is an exponent 
that depends on $q,n$ and $\tau$, according to \hyperref[lem:exponent]{Lemma~\ref*{lem:exponent}}. 
The notation $q^{i - 1} \cdot L_i(\tx)$ denotes a disjunctive repetition of the literal $L_i(\tx)$ for $q^{i - 1}$ times.

\medskip

The idea of this construction is similar to the one for the Viterbi semiring.
While $\epsilon$-characteristic sentences work for $\Vit$-interpretations
where the differences of two distinct values are at least $\epsilon$,
$q$-characteristic sentences work for $\Nat$-interpretations with values less than $q$.
With this in mind, we can explain the construction of $\psi_q(\tx)$ as follows.
If all values in $\pi_A$ are less than $q$, we can picture the value
\[\pi_A \intpr{\psi_q(\ta)} = \sum_{i = 1}^k q^{i - 1} \pi_A(L_i(\ta)),\]
to be in a number system with radix $q$, hence the values $\pi_A(L_i(\ta))$ can be seen as digits.
Thus, it is immediately clear that for any $\Nat$-interpretation $\pi_B$
with universe enumerated by $\tb$ and values less than $q$, 
$\pi_A \intpr{\psi_q(\ta)} = \pi_B \intpr{\psi_q(\tb)}$ implies that $\ta \isoto \tb$ is an isomorphism
between $\pi_A$ and $\pi_B$, since the corresponding ``digits''
$\pi_A(L_i(\ta))$ and $\pi_B(L_i(\tb))$ for each $1 \le i \le k$ have to be the same.

The only remaining problem is the fact that existential quantifiers
$\exists x_1 \ldots \exists x_n$ from $\chi_{\pi_A, q}$ are interpreted as a sum in $\Nat$.
Thus, the value $\pi_A \intpr{\chi_{\pi_A, q}}$ is not induced by a single variable assignment $\ta$.
The exponent $e$ is used to separate the contributions of different variable assignments to the sum
on the basis of \hyperref[lem:exponent]{Lemma~\ref*{lem:exponent}}. 

\begin{theorem}
\label{thm:qcharacteristic}
Let $\pi_A$ and $\pi_B$ are finite $\Nat$-interpretations with values less than $q$. 
Then $\pi_A \intpr{\chi_{\pi_A, q}} = \pi_B \intpr{\chi_{\pi_A, q}}$ implies that $\pi_A \iso \pi_B$.
\end{theorem}

\begin{proof}
Clearly, $\phi(\tx)$ takes care of the number of elements, hence we can assume $\ta$ and $\tb$ 
enumerate the universes of $\pi_A$ and $\pi_B$. Now, we have
\[\pi_A \intpr{\chi_{\pi_A, q}} = \sum_{\sigma \in S_n} \pi_A \intpr{\psi_q(\sigma(\ta))}^e = \sum_{\sigma \in S_n} \pi_B \intpr{\psi_q(\sigma(\tb))}^e = \pi_B \intpr{\chi_{\pi_A, q}}.\]
Recall that $\psi_q(\tx)$ is constructed as a number with $k$ ``digits'', where the digits are the values of the literals 
$\pi_A(L_i(\ta))$ and $\pi_B(L_i(\tb))$, which are bounded by $q$.
Hence, $\pi_A \intpr{\psi_q(\sigma(\ta))}$ and $\pi_B \intpr{\psi_q(\sigma(\tb))}$
are less than $c\coloneqq q^k$, which only depends on $q$, $n$ and $\tau$. 
By \hyperref[lem:exponent]{Lemma~\ref*{lem:exponent}}, there is a sufficiently large $e$ so that 
 $\sum_{\sigma \in S_n} \pi_A \intpr{\psi_q(\sigma(\ta))}^e = \sum_{\sigma \in S_n} \pi_B \intpr{\psi_q(\sigma(\tb))}^e$
implies that both sums share the same summands. In particular, there are permutations
$\sigma_A, \sigma_B \in S_n$ such that $\pi_A \intpr{\psi_c(\sigma_A(\ta))} = \pi_B \intpr{\psi_c(\sigma_B(\tb))}$.

Thanks to the construction of $\psi_c(\tx)$, this yields $\pi_A(L_i(\sigma_A(\ta)) = \pi_B(L_i(\sigma_B(\tb))$ for all 
literals $L_i$ of $\Lit_n(\tau)$. Thus, $\sigma_B \circ (\ta \mapsto \tb) \circ \sigma_A^{-1}$ is an isomorphism from $\pi_A$ to $\pi_B$.
\end{proof}

We can further use the $q$-characteristic sentences also for $\Nat[X]$-interpretations instead of $\Nat$-interpretations. 
Let $X_k = \{x_1, \ldots, x_k\}$, and let $\Nat[X_k](C, n)$ denote the set of polynomials $p \in \Nat[X_k]$ 
with coefficients less than $C$ and exponents less than $n$.
If we choose a suitable variable assignment $X_k \to \Nat$, we can obtain a homomorphism
that assigns unique values to all polynomials in $\Nat[X_k](C, n)$.

\begin{lemma}
The variable assignment $x_i \mapsto C^{n^{i - 1}}$ for $1 \le i \le k$
defines a homomorphism $h \from \Nat[X_k] \to \Nat$ which induces a bijection
from $\Nat[X_k](C, n)$ to $\{0, \ldots, c - 1\} \subseteq \Nat$ where $c \coloneqq C^{n^k}$.
\end{lemma}

\begin{proof}
We proceed by induction on the number of variables $k \in \Nat$.
The base case $k = 0$ is trivial, since $\Nat[\emptyset] \iso \Nat$ and
the empty assignment induces the corresponding isomorphism.
For $k > 0$, notice that $\Nat[X_k] \iso \Nat[X_{k-1}][x_k]$. Hence, each $p \in \Nat[X_k](C, n)$ may be written as
\[p = \sum_{i = 0}^{n - 1} q_i x_k^i, \quad\text{where } q_i \in \Nat[X_{k-1}](C, n).\]
Thus, applying the induced homomorphism $h$ yields
\[h(p) = \sum_{i = 0}^{n - 1} h(q_i) h(x_k)^i = \sum_{i = 0}^{n - 1} h'(q_i) (C^{n^{k - 1}})^i,\]
where $h' \from \Nat[X_{k-1}] \to \Nat$ is induced by $x_i \mapsto C^{n^{i - 1}}$ for $1 \le i < k$.
By induction hypothesis the restriction $h'\vert_{\Nat[X_{k-1}](C, n)}$ is a bijection
from $\Nat[X_{k-1}](C, n)$ to $\{0, \ldots, C^{n^{k - 1}} - 1\}$.
Clearly, $h(p)$ may be seen as a number with $n$ digits $h'(q_i) \in \{0, \ldots, C^{n^{k - 1}} - 1\}$ for $0 \le i < n$
in the number system with radix $C^{n^{k - 1}}$. Thus, any number in $\{0, \ldots, C^{n^k} - 1\}$
can be uniquely represented as $h(p)$ for $p \in \Nat[X_k](C, n)$, which completes the proof.
\end{proof}

\begin{corollary}
For finite $\Nat[X]$-interpretations $\pi_A$ and $\pi_B$
whose values are contained in $\Nat[X](C, n)$,
$\pi_A \intpr{\chi_{\pi_A, c}} = \pi_B \intpr{\chi_{\pi_A, c}}$ implies
$\pi_A \iso \pi_B$ with $c \coloneqq C^{n^{\vert X \vert}}$.
\end{corollary}

\begin{proof}
Transform $\pi_A$ and $\pi_B$ to $\Nat$-interpretations
$\pi_A' \coloneqq h \circ \pi_A$ and $\pi_B' \coloneqq h \circ \pi_B$
by applying the homomorphism from above.
\hyperref[lem:fundamental]{The fundamental property} yields
$\pi_A' \intpr{\chi_{\pi_A, c}} = \pi_B' \intpr{\chi_{\pi_A, c}}$.
Since $h\vert_{\Nat[X](C, n)}$ is a bijection from $\Nat[X](C, n)$ to $\{0, \ldots, c\}$,
the values of $\pi_A'$ and $\pi_B'$ are less than $c$, hence we can invoke
\hyperref[thm:qcharacteristic]{Theorem~\ref*{thm:qcharacteristic}} to conclude
$\pi_A' \iso \pi_B'$. Now, the injectivity of $h$ on $\Nat[X](C, n)$ yields $\pi_A \iso \pi_B$.
\end{proof}

Similarly to the implications of \hyperref[prop:epscharacteristic]{Proposition~\ref*{prop:epscharacteristic}} on the Viterbi semiring $\Vit$,
we conclude that finite $\Nat[X]$-interpretations $\pi_A$ are characterised by
a the set $X_{\pi_A} \coloneqq \{\chi_{\pi_A, c} \mid c \in \Nat, c > 1\}$ of characteristic sentences.
The obvious consequence is that no counterexamples exist on $\Nat[X]$.

\begin{theorem}
For finite $\Nat[X]$-interpretations $\pi_A$ and $\pi_B$, $\pi_A \equiv \pi_B$ implies $\pi_A \iso \pi_B$.
\end{theorem}

Note that the results from this section provide an insight into
the properties of semirings that facilitate the construction of characteristic sentences.
We shall use these observations to separate elementary equivalence from isomorphism  
in semirings that break those properties.

\section{Cancellation}

One of the crucial properties for the characteristic sentences to work is cancellation.
We observe that $\Nat[X]$ and $\Nat$ allow additive and multiplicative cancellation and
$\Vit$ allows multiplicative cancellation, that is, $ab = ac$ implies $b = c$ for all $a \neq 0$.
Looking at $\psi_{\epsilon}(\tx)$, we notice that cancellation is important to ensure that
the contributions of all literals are preserved in a conjunction.
Our analysis will show that there are non-isomorphic but elementary equivalent interpretations 
for a large class of semirings that break cancellation.

\begin{definition}
\label{def:cancellation}
Let $K$ be an idempotent semiring.
A witness that $K$ breaks cancellation is a triple
$a, b, c \in K \setminus \{0\}$ such that
\begin{alignat*}{1}
(1)\quad&a + b = a + c = a \quad\text{and} \\
(2)\quad&ab = ac, \text{ but } b \neq c. \footnotemark 
\end{alignat*}
\end{definition}

For any such triple, we define the following two non-isomorphic $K$-interpretations.
\footnotetext{Strictly speaking, condition (2) suffices for breaking cancellation.
Condition (1) imposes a further condition on the witness, which will be needed below in
the proof.} 

\vspace{\baselineskip}
\begin{minipage}{\linewidth}
\centering
$\pi_b:\quad$
\begin{tabular}{c | c | c |}
$A$ & $R$ & $\neg R$ \\ \hline
$d$ & $a$ & $0$ \\
$e$ & $b$ & $0$ \\
\end{tabular}
$\quad\quad\quad\pi_c:\quad$
\begin{tabular}{c | c | c | c | c |}
$A$ & $R$ & $\neg R$ \\ \hline
$d$ & $a$ & $0$ \\
$e$ & $c$ & $0$ \\
\end{tabular}
\end{minipage}
\vspace{\baselineskip}

\begin{lemma}
\label{lem:cancellation}
The $K$-interpretations $\pi_b$ and $\pi_c$ are elementarily equivalent.
\end{lemma}

\begin{proof}
Consider the $\Bool[X]$ interpretation

\vspace{\baselineskip}
\begin{minipage}{\linewidth}
\centering
$\pi:\quad$
\begin{tabular}{c | c | c |}
$A$ & $R$ & $\neg R$ \\ \hline
$d$ & $x$ & $0$ \\
$e$ & $y$ & $0$ \\
\end{tabular}
\end{minipage}
\vspace{\baselineskip}

Let $h_b, h_c \from \Bool[X] \to K$ be the unique homomorphisms induced
by $x \mapsto a, y \mapsto b$ and $x \mapsto a, y \mapsto c$ respectively.
Obviously, $\pi_b = h_b \circ \pi$ and $\pi_c = h_c \circ \pi$, hence, for each sentence $\psi \in \FO(\{R\})$,
\hyperref[lem:fundamental]{the fundamental property} yields
$\pi_b \intpr{\psi} = h_b(\pi \intpr{\psi})$ and $\pi_c \intpr{\psi} = h_c(\pi \intpr{\psi})$.
In fact, if we set $p \coloneqq \pi \intpr{\psi}$, we have $\pi_b \intpr{\psi} = h_b(p)$ and $\pi_c \intpr{\psi} = h_c(p)$,
hence both interpretations evaluate the same polynomial $p$ under their own homomorphism.

It remains to show that $h_b(p) = h_c(p)$. 
The automorphism $\overline{h}$ of $\Bool[X]$ induced by swapping the variables $x$ and $y$
yields the $\Bool[X]$-interpretation

\vspace{\baselineskip}
\begin{minipage}{\linewidth}
\centering
$\overline{h} \circ \pi = \overline{\pi}:\quad$
\begin{tabular}{c | c | c |}
$A$ & $R$ & $\neg R$ \\ \hline
$d$ & $y$ & $0$ \\
$e$ & $x$ & $0$ \\
\end{tabular}
\end{minipage}
\vspace{\baselineskip}

Clearly, $\pi \iso \overline{\pi}$ by swapping $d$ and $e$, hence
$p = \pi \intpr{\psi} = \overline{\pi} \intpr{\psi} = \overline{h}(\pi \intpr{\psi}) = \overline{h}(p)$.
In other words, $p$ is invariant under swapping variables, so for each
pair $i,j$ we have that $x^iy^j \in p$ if, and only if, $x^jy^i \in p$.
In particular, $x^i \in p \Leftrightarrow y^i \in p$ $(*)$.

Since $p$ is finite and all exponents are less than some $d \in \Nat$, we may write $p$ as
\[p = \sum_{i, j < d, x^iy^j \in p} x^iy^j = \sum_{i < d, x^i \in p} x^i + \sum_{j < d, y^j \in p} y^j + \sum_{0 < i, j < d, x^iy^j \in p} x^iy^j.\]
Set $I \coloneqq \{i < d \mid x^i \in p\} \stackrel{*}{=} \{j < d \mid y^j \in p\}$ and
$M \coloneqq \{x^iy^j \mid 0 < i, j < d, x^iy^j \in p\}$. Then,
\[p = \sum_{i \in I} (x^i + y^i) + \sum_{m \in M} m.\]
For each $m \in M$, $h_b(m) = a^ib^j = a^ic^j = h_c(m)$ due to condition (2) and $i > 0$.
More precisely, if $i > 0$, we can invoke (2) inductively to transform $ab^j$ into $ac^j$ due to commutativity of multiplication. We now invoke condition (1) for each $i \in I$. For $z \in \{b, c\}$ and using idempotence of $K$ (i), this yields
\begin{alignat*}{1}
a^i + z^i &\stackrel{(1)}{=} (a + z)^i + z^i = \sum_{j = 0}^i a^{i - j}z^j + z^i = \sum_{j = 0}^{i - 1} a^{i - j}z^j + z^i + z^i \\
&\stackrel{(i)}{=} \sum_{j = 0}^{i - 1} a^{i - j}z^j + z^i = a^i.
\end{alignat*}
Hence, $h_b(x^i + y^i) = a^i + b^i = a^i = a^i + c^i = h_c(x^i + y^i)$ for each $i \in I$ holds as well. Together, we have
\[h_b(p) = \sum_{i \in I} h_b(x^i + y^i) + \sum_{m \in M} h_b(m) = \sum_{i \in I} h_c(x^i + y^i) + \sum_{m \in M} h_c(m) = h_c(p),\]
which completes the proof, since $\psi$ was arbitrary and $\pi_b \intpr{\psi} = h_b(p) = h_c(p) = \pi_c \intpr{\psi}$.
\end{proof}

\hyperref[lem:cancellation]{Lemma~\ref*{lem:cancellation}} can be applied to many important semirings, 
such as $\Bool[X]$ itself and $\Sorb[X]$ with an appropriate choice of $a$, $b$ and $c$.

\begin{theorem}
For $X \supseteq \{x, y\}$, there exists a pair of elementarily equivalent,
but non-isomorphic $K$-interpretations in the shape of $\pi_b$ and $\pi_c$
for both $K = \Bool[X]$ and $K = \Sorb[X]$.
\end{theorem}

\begin{proof}
For $\Bool[X]$, choose
$a \coloneqq x + y + x^2 + xy + y^2$, $b \coloneqq x^2 + y^2$ and $c \coloneqq x^2 + xy + y^2$ 
to obtain the following pair of $\Bool[X]$-interpretations.

\vspace{\baselineskip}
\begin{minipage}{\linewidth}
\centering
$\pi_b:\quad$
\begin{tabular}{c | c | c |}
$A$ & $R$ & $\neg R$ \\ \hline
$d$ & $x + y + x^2 + xy + y^2$ & $0$ \\
$e$ & $x^2 + y^2$ & $0$ \\
\end{tabular}
$\quad\quad\pi_c:\quad$
\begin{tabular}{c | c | c | c | c |}
$A$ & $R$ & $\neg R$ \\ \hline
$d$ & $x + y + x^2 + xy + y^2$ & $0$ \\
$e$ & $x^2 + xy + y^2$ & $0$ \\
\end{tabular}
\end{minipage}
\vspace{\baselineskip}

To prove the desired properties, we only have to check conditions (1) and (2)
from \hyperref[def:cancellation]{Definition~\ref*{def:cancellation}}
and then invoke \hyperref[lem:cancellation]{Lemma~\ref*{lem:cancellation}}.
Condition (1) is obvious, and for (2), it suffices to expand the products $ab$ and $ac$ to calculate that
\[ab = ac = x^3 + xy^2 + x^2y + y^3 + x^4 + x^3y + x^2y^2 + xy^3 + y^4.\]

An analogous construction works for the semiring of absorptive polynomials $\Sorb[X]$.
Choose $a \coloneqq x + y$, $b \coloneqq x^2 + y^2$ and $c \coloneqq x^2 + xy + y^2$.
Note that this is the same as the counterexample on $\Bool[X]$ above after applying absorption,
which yields the pair of $\Sorb[X]$-interpretations depicted below.

\vspace{\baselineskip}
\begin{minipage}{\linewidth}
\centering
$\pi_b:\quad$
\begin{tabular}{c | c | c |}
$A$ & $R$ & $\neg R$ \\ \hline
$d$ & $x + y$ & $0$ \\
$e$ & $x^2 + y^2$ & $0$ \\
\end{tabular}
$\quad\quad\pi_c:\quad$
\begin{tabular}{c | c | c | c | c |}
$A$ & $R$ & $\neg R$ \\ \hline
$d$ & $x + y$ & $0$ \\
$e$ & $x^2 + xy + y^2$ & $0$ \\
\end{tabular}
\end{minipage}
\vspace{\baselineskip}

It is easy to verify that $a$, $b$ and $c$ also satisfy (1) and (2), which completes the proof.
\end{proof}

\section{Conclusion and outlook}

Our analysis of first-order axiomatisations and elementary equivalence 
of finite semiring interpretations has revealed some remarkable differences
between semiring semantics and classical Boolean semantics.
Depending on the underlying semiring, there may
exist finite semiring interpretations that are elementarily equivalent without being
isomorphic. Indeed, this phenomenon happens already in very simple cases such as
for min-max semirings with three elements. On the other side, there are
relevant semirings, used for instance in provenance analysis in databases such
as the tropical semiring or the Viterbi semiring, where every finite interpretation
is first-order axiomatisable, and in fact even by a finite set of axioms. However,
and this is again an interesting difference to Boolean semantics, a finite
axiomatisation does not imply an axiomatisation by a single axiom.

Also for the semirings of polynomials, fundamental for a general provenance
analysis that reveals which combinations of atomic facts are responsible
for the truth of a logical statement, the picture is not unique. While the 
most general semiring $\N[X]$, freely generated by $X$, admits 
axiomatisations of all finite interpretations, so that elementary equivalence
implies isomorphism, this is not the case for the semirings 
$\PosBool[X]$, $\Sorb[X]$, $\Bool[X]$ and $\Why[X]$ which are 
universal for important subclasses of semirings.

In the study of elementary equivalence for semiring semantics,
it turns out that there is no \emph{straightforward} adaptation of 
\EF{} games, or their generalisations such as
Hellas bijective pebble game \cite{Hella92}, to semiring interpretations.
Whatever the specific protocol of possible moves in such games may be, they
always result in a \emph{localisation}, in the sense that some tuples
in the two structures are picked that are indistinguishable on the atomic level.
As shown by the very simple example of the interpretations $\pi_{PQ}$ and $\pi_{QP}$
at the beginning of \hyperref[sec:separation]{Sect.~\ref*{sec:separation}}, this is not possible
in semiring semantics. Although the two interpretations are elementarily 
equivalent, no element of the first ``looks the same'' as any element of the second,
so any kind of localisation would result in a winning play of the Spoiler.
It is an intriguing open question how elementary equivalence of semiring interpretations
can be captured by a different notion of comparison games or back-and-forth systems à la \Fraisse.
This not being available (yet), we have established elementary equivalence by
different methods, based on homomorphisms, which we believe to be of independent interest.
 
There are many other model-theoretic issues that deserve to be studied in semiring semantics.
While we have limited ourselves here to \emph{finite} semiring interpretations, 
the study of semiring semantics over infinite universes is of course very interesting as well.
It requires certain restrictions on the underlying semirings, concerning existence and
appropriate algebraic properties of infinite sums and products, but there are useful
semirings satisfying such properties. A particularly interesting question is what kind of
compactness and preservation results are possible in such contexts.

We finally remark that an altogether different approach to semiring interpretations would
consider them as two-sorted structures, one sort being a finite or infinite structure
(or just a set), the second one consisting of the semiring, with functions
from the first to the second sort giving the semiring interpretation of the literals.
This is very similar to the approach of metafinite model theory \cite{GraedelGur98},
and to get a reasonable logical theory it is important that the elements of
the second sort, here the semiring, are treated differently than the elements
of the first sort. In particular, quantifiers should range only over the first sort,
and operations on the second sort are just the algebraic semiring operations
and their aggregates, together with equalities between terms. Such an
approach is certainly useful for a number of questions and permits the study of semiring 
interpretations via classical Boolean semantics. In particular, once semiring values
are directly accessible in the logic, the construction of characteristic 
sentences, axiomatising a finite structure up to isomorphism, can easily be translated into
such a setting. However, such an internalisation of the semiring values,
from the meta-level of truth values into the structures under consideration,
does not really capture the essence of semiring semantics. 



\end{document}